\documentclass[10pt]{amsart}
\usepackage{amsbsy}
\usepackage{graphicx,epsfig,subfigure,psfrag}
\begin{document}

\newtheorem{theorem}{Theorem}
\newtheorem{proposition}{Proposition}
\newtheorem{lemma}{Lemma}
\newtheorem{corollary}{Corollary}
\newtheorem{definition}{Definition}
\newtheorem{remark}{Remark}
\newcommand{\tex}{\textstyle}
\numberwithin{equation}{section} \numberwithin{theorem}{section}
\numberwithin{proposition}{section} \numberwithin{lemma}{section}
\numberwithin{corollary}{section}
\numberwithin{definition}{section} \numberwithin{remark}{section}
\newcommand{\ren}{\mathbb{R}^N}
\newcommand{\re}{\mathbb{R}}
\newcommand{\n}{\nabla}
\newcommand{\p}{\partial}
\newcommand{\iy}{\infty}
\newcommand{\pa}{\partial}
\newcommand{\fp}{\noindent}
\newcommand{\ms}{\medskip\vskip-.1cm}
\newcommand{\mpb}{\medskip}
\newcommand{\AAA}{{\bf A}}
\newcommand{\BB}{{\bf B}}
\newcommand{\CC}{{\bf C}}
\newcommand{\DD}{{\bf D}}
\newcommand{\EE}{{\bf E}}
\newcommand{\FF}{{\bf F}}
\newcommand{\GG}{{\bf G}}
\newcommand{\oo}{{\mathbf \omega}}
\newcommand{\Am}{{\bf A}_{2m}}
\newcommand{\CCC}{{\mathbf  C}}
\newcommand{\II}{{\mathrm{Im}}\,}
\newcommand{\RR}{{\mathrm{Re}}\,}
\newcommand{\eee}{{\mathrm  e}}
\newcommand{\LL}{L^2_\rho(\ren)}
\newcommand{\LLL}{L^2_{\rho^*}(\ren)}
\renewcommand{\a}{\alpha}
\renewcommand{\b}{\beta}
\newcommand{\g}{\gamma}
\newcommand{\G}{\Gamma}
\renewcommand{\d}{\delta}
\newcommand{\D}{\Delta}
\newcommand{\e}{\varepsilon}
\newcommand{\var}{\varphi}
\newcommand{\lll}{\l}
\renewcommand{\l}{\lambda}
\renewcommand{\o}{\omega}
\renewcommand{\O}{\Omega}
\newcommand{\s}{\sigma}
\renewcommand{\t}{\tau}
\renewcommand{\th}{\theta}
\newcommand{\z}{\zeta}
\newcommand{\wx}{\widetilde x}
\newcommand{\wt}{\widetilde t}
\newcommand{\noi}{\noindent}
\newcommand{\uu}{{\bf u}}
\newcommand{\xx}{{\bf x}}
\newcommand{\yy}{{\bf y}}
\newcommand{\zz}{{\bf z}}
\newcommand{\aaa}{{\bf a}}
\newcommand{\cc}{{\bf c}}
\newcommand{\jj}{{\bf j}}
\newcommand{\ggg}{{\bf g}}
\newcommand{\UU}{{\bf U}}
\newcommand{\YY}{{\bf Y}}
\newcommand{\HH}{{\bf H}}
\newcommand{\GGG}{{\bf G}}
\newcommand{\VV}{{\bf V}}
\newcommand{\ww}{{\bf w}}
\newcommand{\vv}{{\bf v}}
\newcommand{\hh}{{\bf h}}
\newcommand{\di}{{\rm div}\,}
\newcommand{\ii}{{\rm i}\,}
\def\I{{\rm Id}}
\newcommand{\inA}{\quad \mbox{in} \quad \ren \times \re_+}
\newcommand{\inB}{\quad \mbox{in} \quad}
\newcommand{\inC}{\quad \mbox{in} \quad \re \times \re_+}
\newcommand{\inD}{\quad \mbox{in} \quad \re}
\newcommand{\forA}{\quad \mbox{for} \quad}
\newcommand{\whereA}{,\quad \mbox{where} \quad}
\newcommand{\asA}{\quad \mbox{as} \quad}
\newcommand{\andA}{\quad \mbox{and} \quad}
\newcommand{\withA}{,\quad \mbox{with} \quad}
\newcommand{\orA}{,\quad \mbox{or} \quad}
\newcommand{\atA}{\quad \mbox{at} \quad}
\newcommand{\onA}{\quad \mbox{on} \quad}
\newcommand{\ef}{\eqref}
\newcommand{\mc}{\mathcal}
\newcommand{\mf}{\mathfrak}

\newcommand{\ssk}{\smallskip}
\newcommand{\LongA}{\quad \Longrightarrow \quad}
\def\com#1{\fbox{\parbox{6in}{\texttt{#1}}}}
\def\N{{\mathbb N}}
\def\A{{\cal A}}
\newcommand{\de}{\,d}
\newcommand{\eps}{\varepsilon}
\newcommand{\be}{\begin{equation}}
\newcommand{\ee}{\end{equation}}
\newcommand{\spt}{{\mbox spt}}
\newcommand{\ind}{{\mbox ind}}
\newcommand{\supp}{{\mbox supp}}
\newcommand{\dip}{\displaystyle}
\newcommand{\prt}{\partial}
\renewcommand{\theequation}{\thesection.\arabic{equation}}
\renewcommand{\baselinestretch}{1.1}
\newcommand{\Dm}{(-\D)^m}

\title
{\bf The ${\bf P}$-Laplace equation in domains \\with multiple crack
section via pencil operators}

\author{Pablo~\'Alvarez-Caudevilla and Victor~A.~Galaktionov}

\address{Universidad Carlos III de Madrid,
Av. Universidad 30, 28911-Legan\'es, Spain -- Work phone number: +34-916249099}
\email{pacaudev@math.uc3m.es}

\address{Department of Mathematical Sciences, University of Bath,
 Bath BA2 7AY, UK}
\email{vag@maths.bath.ac.uk}


\thanks{This work has been partially supported by the Ministry of Economy and Competitiveness of
Spain under research project MTM2012-33258.}

\date{\today}




\begin{abstract}
 The $p$-Laplace equation
  $$
  \n \cdot (|\n u|^{p-2} \n u)=0 \whereA p>2,
  $$
 in a bounded domain $\O \subset \re^2$,
 with inhomogeneous Dirichlet
  conditions on the smooth boundary $\p \O$ is
 considered. In addition, there is a finite collection of  curves
  $$
 \Gamma = \Gamma_1\cup...\cup\Gamma_m \subset \O, \quad
 \mbox{on which we assume homogeneous Dirichlet conditions} \quad  u=0,
  $$
  modeling a multiple crack formation, focusing at the origin $0 \in \O$.
  This  makes the above quasilinear elliptic problem overdetermined. Possible
  types of the behaviour of solution $u(x,y)$ at the tip $0$ of such
  admissible multiple
  cracks, being a ``singularity" point, are described, on the basis of
   blow-up
  scaling techniques and
 a ``nonlinear eigenvalue problem".
  Typical types of admissible cracks are
  shown to be governed by nodal sets of a countable family of
  {\em
nonlinear eigenfunctions}, which are obtained
via
   branching from harmonic polynomials that occur for  $p=2$.
Using a combination of analytic and numerical methods, saddle-node
bifurcations in $p$ are shown to occur for those nonlinear
eigenvalues/eigenfunctions.

\end{abstract}

\maketitle

\section{Introduction}
 \label{S1}

\subsection{Models and preliminaries}


\noindent We study
 solutions of the \emph{$p$-Laplace equation} with Dirichlet
boundary conditions in a bounded smooth domain $\O \subset \re^2$
\begin{equation}
\label{pa}
   \left\{\begin{array}{cc} \n \cdot (|\n u|^n \n u)  =0 &
\hbox{in}\quad \O,\\
    u=f(x,y) &  \hbox{on}\quad \G,\\
   u=g(x,y) &  \hbox{on}\quad \p\O,\\
\end{array} \right.
\end{equation}
where $n>0$ is a fixed
exponent ($n=p-2$,  for the standard $p$-Laplacian but we use the parameter $n$ for convenience in our subsequent branching analysis) 
and $f(x,y)$ and $g(x,y)$
are  given smooth functions.
In our particular case, $\O$
is assumed to have  a multiple crack $\G$, as a finite collection
of $m \ge 1$ curves
 \be
 \label{curv1}
  \G= \G_1 \cup \G_2 \cup...\cup \G_m \subset \O  \quad \mbox{such that each $\G_j$ passes through the origin $0 \in \O$ only}.
   \ee
The origin is then  the tip of this crack. 

\begin{figure}[htp]

\includegraphics[scale=0.4]{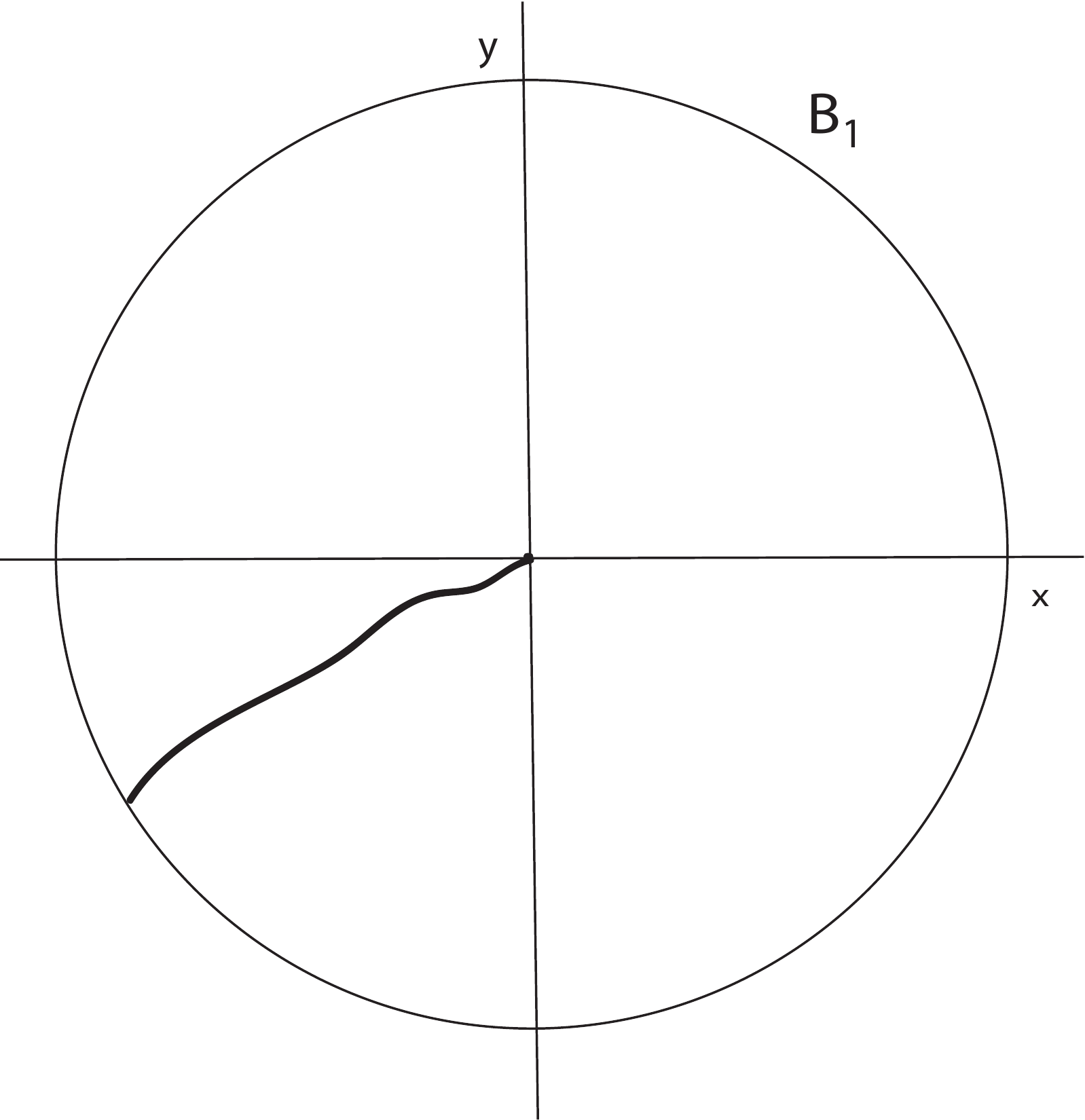}

\vskip -0.2cm \caption{ \small A formal one-crack model, $m=1$.}
\label{figcrack}
\end{figure}

Moreover, we assume that, near the origin, in the
lower half-plane $\{y<0\}$, all cracks asymptotically  take a
straight line form, i.e., as shown in Figure \ref{figcrack},
 \be
 \label{str1}
 \G_k:  \,\, x=\a_k (-y)(1+o(1)), \,\, y \to 0, \,\,\,
 k=1,2,...,m,\,\,
 \mbox{where} \,\,\,
 \a_1<\a_2<...< \a_m
  \ee
  are given constants.
Thus, the precise statement of the problem assumes that \ef{str1}
describes {\em all the admissible cracks near the origin}, i.e.,
{\em no other straight-line cracks are considered}. Indeed, through our analysis it will be determined that these are the only 
type of admissible cracks.

In our basic model, without loss of any generality (and
to simplify the analysis), we assume a homogeneous Dirichlet
condition at the crack:
\be
 \label{Dir1}
 u=0 \onA \G,
  \ee
  that makes the problem overdetermined. 
  So that, not {\em any} type of such
   multiple cracks \ef{curv1}, \ef{str1} are admissible.
   
 For this problem, we ascertain important qualitative information about the behaviour
  of the solutions of problem \eqref{pa}. 
 Especially important is the analysis close to the tip of the crack $\G$, at least when the parameter $n$ is very close to zero (or $p$ very close to 2).  
 
 In particular the behaviour 
 of the solutions at the singularity boundary points are described by some blow-up scaling techniques but also 
 by a nonlinear eigenvalue problem through a branching argument.
 In other words we use an appropriate change of variable to transform the equations involved into the so-called 
 pencil operators which reduces the problem to solve a 1$D$ spectral problem, allowing us to 
 ascertain such an asymptotical behaviour at the tip of the crack.

Thus, using Pencil Operator Theory we first analyse the problem when $n=0$ reducing \eqref{pa} to the well known Laplace 
equation. Most of the results presented for this problem are already known; such us the expression of its harmonic polynomials
and their properties. However, we find convenient to show these properties within the pencil operator framework from which this work has been developed. 
Something that we believe is not widely known.

Furthermore, this analysis will be very important in the subsequent regarding the results of the $p$-Laplace problem 
\eqref{pa}. This problem is nonlinear and we cannot apply directly the ideas used for the linear case when $n=0$. 
Therefore, in order to arrive at the results we perform a branching analysis from the 
solutions of the problem when $n=0$ for which we have plenty of information.

\ssk

\subsection{Laplace equation and pencil operators}


Thus, first for $n=0$ we describe the behaviour of the solutions for the Laplace problem 
 \be
 \label{Lap1}
 \D u=0 \inB \,\,\O\;(=B_1), \quad u=f \,\,(\not \equiv 0) \onA \p
 \O, \quad \mbox{with  (\ref{Dir1})}.
  \ee
at the tip of the crack (normally the origin), being a singularity boundary point. 

To do so, blow-up scaling techniques 
and spectral theory of pencils of non self-adjoint operators are used. 
In particular, we perform a rescaling/change of variable 
of the form
  \be
  \label{scal}
  \tex{  z= x/(-y) \quad \hbox{and}\quad  \t = - \ln
(-y) \forA y<0. } \ee
 This rescaling corresponds to a blow-up
scaling near the origin. In fact, this blow-up analysis assumes a
kind of elliptic ``evolution" approach for elliptic problems.
We actually move the singularity point at the origin into an
asymptotic convergence when $\tau \to \infty$.

Applying this rescaling (and then the separation variables method) 
we transform the Laplace equation into a pencil of non-self-adjoint operator, in particular for this 
case, a quadratic pencil operator of the form
$$(1+z^2)(\psi^*)''+2(\l+1)z (\psi^*)' +\l(\l+1)
 \psi^*=0,
$$
for which we obtain two families of eigenfunctions 
\be
\label{harm}
\{\psi_{k_+}^*(z)\equiv \psi_{k,1}^*(z)\} \quad \hbox{and}\quad \{\psi_{k_-}^*(z)\equiv  \psi_{k-1,2}^*(z)\}\quad \hbox{(defined in Section\;\ref{SLap})},
\ee
associated with two families of eigenvalues
\be
\label{eifam}
 \l_k^+=-l, \quad l=1,2,3,... \andA \l_k^-=-l-1, \quad
  l=0,1,2,3,... \,,
  \ee
Consequently, the application of  the theory of pencil operators allows us to analyse the linear equation
\eqref{Lap1}, with the geometrical structure assumed in this work.

Pencil operators are  polynomials of the form $$ A(\l):=A_0+\l
A_1+\cdots+\l^n A_n,$$
 where $\l \in
{\mathbb C}$ is a spectral parameter and the coefficients $A_i$,
with $i=0,1,\cdots, n$, are linear operators acting on a Hilbert
space. Pencil operator theory appeared and was crucially used in
the regularity and asymptotic analysis of elliptic problems in a
seminal paper by Kondrat'ev \cite{Kon2} and also for parabolic
problems in \cite{Kon1}, where spectral problems, that are
nonlinear (polynomial) in the spectral parameter $\l$, occurred.
Later on, Mark Krein and Heinz Langer \cite{KL} made a fundamental
contribution to this theory analyzing the spectral theory for
strongly damped quadratic operator pencils.

Thus, using this well-developed spectral
theory of non-self-adjoint operators, we are
able to ascertain information about the solutions at the tip of
the crack $\G$ for the Laplace problem \eqref{Lap1}. Then, we deal with eigenfunctions \eqref{harm} of a
quadratic pencil of operator
 and not with a standard Sturm--Liouville problem
 (we keep this ``blow-up scaling logic" for the analysis of the $p$-Laplace equation as well).

Consequently, and in particular for this simpler problem, we obtain that all the solutions with cracks at 0 must have the expression
 \be
 \label{decola}
  \tex{
  u(x,y)=w(z,\t)= \sum_{(k \ge l)} \eee^{-k \t} [c_k\psi_{k,1}^*(z)+ d_k \psi_{k-1,2}^*(z)],
   \quad \mbox{with} \quad
c_l^2+d_l^2 \ne 0,
  }
  \ee
  where $\psi_{k,1}^*(z)$ and $\psi_{k-1,2}^*(z)$ are two families of harmonic polynomials
  re-written in terms of a rescaled variable $z$ after performing the scaling \eqref{scal}.

 Therefore, the main result established for the Laplace problem shows that in the first leading terms while approaching
the origin, a linear combination of two families of eigenfunctions
as classic harmonic polynomials. Then,
 
 \begin{itemize}
 \item If all $\{\a_k\}$ in \ef{str1} do not coincide with
all $m$ subsequent zeros of any nontrivial linear combination
 \be
 \label{dm1}
c_l \psi_{l,1}^*(z)+d_l \psi^*_{l-1,2}(z), \quad \mbox{with} \quad
c_l^2+d_l^2 \ne 0,
 \ee
 then the multiple crack problem \ef{Lap1}
cannot have a solution for any boundary Dirichlet data $f$ on $\p
\O$.
\item However, if that is not the case, i.e. there is some $l$ for which the $m$ zeros of 
\eqref{dm1} coincide with $\{\a_k\}$, then there exists a solution $u$ and
$$|u(x,y)|=O\left(|x,y|^l\right)\quad \hbox{as}\quad (x,y)\to (0,0).$$
\end{itemize}

Moreover and obviously, restricting to $\G$
all types of admissible crack-containing expansions \ef{decola}
(with closure in any appropriate functional space)  {\em fully
describe all types of boundary data, which lead to the desired
crack formation at the origin}.

Thus, typical types of admissible 
cracks are shown to be governed by nodal sets of a countable family of harmonic polynomials, which 
are represented by pencil eigenfunctions, instead of their classical representation via a standard Sturm-Liouville problem;
see Section\;\ref{SLap}.  The analysis carried out for the case $n=0$ (the Laplace equation) 
will show that actually those types of cracks \eqref{str1} are the only
  admissible ones.


\subsection{Our main problem and main new results: $p$-Laplace equation}

 For the
quasilinear problem \ef{pa},  the obtained results for the Laplace crack problem \ef{Lap1}
seem to be crucial. However, we need to
change the strategy and base our analysis on a branching argument
at $n=0$. 

The main reason is due to the fact that for the previous Laplace
equation we are able to obtain explicitly two families of
negative eigenvalues associated, respectively, to two families of
eigenfunctions. 
Indeed, from the expressions of the two families
of eigenvalues \eqref{eifam} (see Section\,\ref{SLap} below) and using the
 pencil operator theory \cite{Mar}, we can ascertain the
coefficients of every eigenfunction explicitly as well. 

However,
following the same argument for the nonlinear PDE \eqref{pa} it
is not possible, in general, to get the corresponding families of
eigenvalues and, hence, the associated eigenfunctions.

Therefore, using the knowledge we possess about the eigenfunctions of
the quadratic pencil operator coming from the Laplace equation
\eqref{Lap1}, and after performing the blow-up scaling \eqref{scal},
we carry out a branching analysis at $n=0$ obtaining information about
the solutions/nonlinear eigenfunctions, at least, when $n$ is sufficiently close to zero
and with the geometrical framework under consideration in this work.

Namely, we show that admissible crack configurations also
correspond to nodal sets of {\em nonlinear eigenfunctions}, which
we construct by branching at $n=0$ (or $p=2$) from the (re-written
in terms of $z$) harmonic polynomials. 

In particular, after performing the rescaling \eqref{scal} we transform the $p$-Laplace problem 
\eqref{pa} into a pencil operator (see the details in Section\;\ref{plaplace}) which is indeed a nonlinear eigenvalue problem. Then, 
we construct a set
of uniquely determined nonlinear eigenfunctions $\Psi^*_l(z)$, emanating at
$n=0$ from the harmonic polynomial eigenfunctions $\psi_l^*(z)$ denoted by \eqref{harm} such that
$$ \Psi^*_l(z)= \psi_l^*(z) + n\, \var_l(z) + o(n),$$
and the associated eigenvalues follow
$$\Lambda_ l=\l_l+n \, \mu_l+o(n),$$
where $\var_l(z)$ are unknown functions and $\mu_l$ unknown constants. 

However, even with the previous relationship between the non-linear eigenfunctions $ \Psi^*_l(z)$ and the eigenfunctions $\psi_l^*(z)$
of the quadratic pencil operator associated with Laplace equation \eqref{Lap1} we observe that the ``blow-up" zero structures or nodal sets 
may be very different for the $p$-Laplace equation \eqref{pa}.

Through this local analysis we can only assure
uniqueness when the parameter $n$ is sufficiently close to zero.
Otherwise, the possibility of a countable family of solutions
appears to be very likely. This is summarized by Theorem\;\ref{mainpt}. 
Note that, contrary to what happens for the Laplacian the admissible cracks for the $p$-Laplacian problem \eqref{pa} 
 are governed by nodal sets of a 
 countable family of nonlinear eigenfunctions.

 Therefore, using analytic and numerical methods saddle-node bifurcations are shown to occur for 
 those nonlinear eigenvalues/eigenfunctions, 
 so that there is no existence of real
solutions above a corresponding critical value $n^*$ of the
parameter $n$. Note that a global continuation (when $n$ is away from zero, $n>0$) of such $n$-branches requires a much more 
difficult analytical and even numerical analysis than the one explained here.

Hence, up to that critical value for the parameter $n$ we expect to have a similar nodal set to 
the linear case (the Laplace problem \eqref{Lap1}). However beyond this maximal critical value one can expect 
more complicated nodal sets. Nevertheless, this fact still remains unanswered. 

Furthermore, note that the completeness of the nonlinear eigenfunctions $\{\Psi_l^*(z)\}$, although
it can be expected and it is necessary to complete the classification of the crack configuration for the $p$-Laplacian problem \eqref{pa}, 
is also a very difficult open problem. This is basically due to the complexity expected 
for the nodal sets of those nonlinear eigenfunctions.   

In conclusion, performing a proper rescaling \eqref{scal} in the problem
\eqref{pa}
and using ``nonlinear operator pencil theory", we are able to show
those special linear combinations of  ``nonlinear harmonic
polynomials".
Specifically, we claim that their nodal sets
 play a key role in the general multiple crack problem for various
nonlinear elliptic equations.

\ssk

Though our approach is done in two dimensions, the scaling blow-up
approach applies to $\O$ in $\re^3$ (or any $\ren$), where
spherical ``nonlinear harmonic polynomials" naturally occur so
that their nodal sets (finite combination of nodal surfaces) of
their linear combinations, as above, describe all possible local
structures of cracks concentrating at the origin. However, if $N>2$ 
the possible geometry of the crack is far richer.

Note that domain $\O\setminus \Gamma$ fits in the context of the definition of
smooth cones in $\ren$ (here we assume $N=2$).
 In other words, we say that the crack $\G$ is a smooth cone if it is a set of dimension
$N-1$ in $\re^N$, conical, centered at the origin and
$S^{N-1}\setminus\G$ is a domain with a piecewise $C^2$ boundary;
see \cite{AL} and references therein for any further details.
Moreover, we would just like to mention that the proof therein is based on the
assumption that the imbedding $W^{1,2}( S^{N-1}\setminus\G)$ into
$L^2(S^{N-1}\setminus\G)$ is compact (cf. \cite{Adm}).


\subsection{Further extensions}


For instance, the results obtained for the $p$-Laplacian operator
\eqref{pa} can be extended to the bi-$p$-Laplacian equation
 \be
 \label{d2}
 \D_p^2 u  \equiv \D(|\D|^{n} \D u)=0 \whereA p=2+n,
  \ee
  which leads to much more complicated technical computations, but
  these results are out of the scope of this work. However, we would just like to
  mention that, as performed for the $p$-Laplacian case the corresponding {\em nonlinear
  eigenfunctions} of the nonlinear pencil for \ef{d2} can be
  obtained by branching from harmonic polynomials such as
  eigenfunctions of a {\em polynomial $($quartic$)$ pencil}
of a non self-adjoint operator that occurs for the
  bi-Laplacian in the Dirichlet problem
  \begin{equation}
\label{bila}
\left\{\begin{array}{cc}
    \Delta^2 u =0   &  \hbox{in}\quad \O,\\
    u=f(x,y) &  \hbox{on}\quad \G,\\
    u=g(x,y),\;\; \frac{\p u}{\p {\bf n}}=h(x,y) &  \hbox{on}\quad \p\O,\\
\end{array} \right.
\end{equation}
with the same zero-condition \ef{Dir1} on the multiple cracks. See
\cite{CGbiLap} for any further details, applications, proofs, and
 discussions about this particular problem \eqref{bila}.

\section{The linear case $n=0$: crack distribution via nodal sets of transformed harmonic
 polynomials}
 \label{SLap}

In this section we assume that $n=0$ which leads us  to the multiple
crack Laplace problem \ef{Lap1}.
We now present several results to be used later in the analysis of the $p$-Laplace equation \eqref{pa} 
via a branching analysis with $n>0$.

\subsection{Blow-up scaling and rescaled equation}

First we show the required transformations with which we will obtain the pencil operators 
that will eventually provide us with the behaviour of the solutions at the tip of the crack for the Laplace problem \eqref{Lap1}. 

Thus, assuming the crack configuration as in Figure \ref{figcrack}, we
introduce the following rescaled variables, corresponding to a
``blow-up" scaling near the origin $0$:
  \be
  \label{resc1}
  \mbox{$
u(x,y)= w(z ,\t), \quad \mbox{with} \quad z=  x/{(-y)}
 \quad \mbox{and} \quad  \t = - \ln
(-y) \forA y<0,
 $}
  \ee
   to get the rescaled operator
 \be
 \label{www.20}
\D_{(x,y)}u= \eee^{2 \t} \,\big[ D_\t^2+ D_\t+ 2 z D^2_{z \t} +
(1+z^2)D_z^2 +2z D_z\big]w \equiv \D_{(z ,\t)}w.
 \ee
 
\vspace{0.2cm}

Consequently, from equation \eqref{www.20} and the scaling \eqref{resc1} in a neighbourhood of the origin, we arrive at the
equation
 \be
 \label{www.2}
  w_{\t\t} + w_\t + 2z w_{z  \t} = {\bf A}^* w,\quad \hbox{such that}\quad  {\bf A}^* w \equiv -(1+z
^2) w_{zz }- 2 z w_z .
 \ee
 
 \vspace{0.2cm}
 
 \noindent{\bf Remark.}  We observe that the operator $\AAA^*$ is symmetric in the standard (dual) metric of
$L^2(\re)$,
 \be
 \label{aa1}
 \AAA^* \equiv - D_z[(1+z^2)D_z],
  \ee
  though we are not going to use this. Indeed, for our crack
  purposes, we do not need eigenfunctions of the ``adjoint"
pencil, 
since we are not going to use
eigenfunction expansions of solutions of the PDE \ef{www.2}, where
bi-orthogonal basis could naturally be wanted.

 \vspace{0.2cm}

This blow-up analysis of \ef{www.2}
assumes a kind of  ``elliptic evolution" approach for elliptic
problems, which  is not well-posed in the so-called Hadamard's sense (see \cite{Had} for more details) but,
indeed, can trace out the behaviour  of necessary global orbits
that reach and eventually decay to the singularity point
$(z,\t)=(0,+\iy)$. Then, 
 by the crack condition \ef{Dir1}, we
look for vanishing solutions: in the mean and  uniformly on
compact subsets in $z$,
 \be
  \label{zer1}
   w(z,\tau) \to 0 \asA \t \to +\iy .
    \ee
Therefore, under the rescaling \eqref{resc1}, we have converted
the singularity point at $(0,0)$ into an asymptotic convergence
when $\tau \to \infty$.

Hence,  we are forced to describe a very thin family of
solutions for which
we will describe their possible nodal
sets to settle the multiple crack condition in \ef{Lap1}. This
corresponds to Kondratiev's ``evolution" approach \cite{Kon1,
Kon2} of 1966, though it was there directed to different {\em
boundary point regularity} (and asymptotic expansions) questions,
while the current crack problem assumes studying the behaviour at
an {\em internal} point $0 \in \O$ such as the tip of the multiple
crack under consideration. We will show first that this {\em
internal crack problem} requires polynomial eigenfunctions of
different pencils of linear operators, which were not under
scrutiny in Kondratiev's.

\subsection{Quadratic pencil and its polynomial eigenfunctions}

We now obtain the quadratic pencil operator associated with equation \eqref{www.2}. Remember 
we have arrived at that equation performing the rescaling \eqref{resc1}. Once this pencil operator 
is obtained we will show several spectral properties that are important in our analysis.

Also, we should say that these properties are very well known for a classical Laplacian. However, since 
we are transforming this linear problem we find convenient to include them for the corresponding quadratic 
pencil operator.

 As usual in linear PDE theory, looking for solutions of
\ef{www.2} in separate variables
 \be
\label{YYY.112S}
 w(z , \t)= {\mathrm e}^ {\l \t} \psi^*(z ) \whereA {\rm Re} \, \l
 <0 \,\,\, \mbox{by (\ref{zer1})},
 \ee
yields the eigenvalue problem for a {\em quadratic pencil}
of non self-adjoint operators, 
 \be
 \label{Pen.1}
 {\bf B}_\l^* \psi^* \equiv  \{\l(\l+1)I+ 2 \l z  D_z  - {\bf A}^*\} \psi^*=0
 \,\, \mbox{or}\,\, (1+z^2)(\psi^*)''+2(\l+1)z (\psi^*)' +\l(\l+1)
 \psi^*=0.
   \ee
   
 \ssk
 
 \noi{\bf Remark.} The  second-order operator ${\bf A}^*$ denoted by \eqref{aa1} is singular  at the
infinite points $z= \pm \infty$, so this is a singular quadratic
pencil eigenvalue problem. 

Moreover, since the linear first-order operator
in \ef{Pen.1}, $z D_z$, is not symmetric
 in $L^2$, we are not obliged to attach the whole operator to any
 particular functional space.
 Therefore, the behaviour as $z \to \infty$  is not that
crucial, and any $L^2_\rho$-space setting  with $\rho(z) \sim
{\mathrm e}^{- a z^2}$ (or ${\mathrm e}^{-a |z|}$), $a>0$ small,
would be enough. 
Indeed, if the solution of the problem \eqref{Lap1} is smooth in
certain weighted spaces $H^1_\rho$ or $L^2_\rho$, we claim that,
by classic ODE theory, then the eigenfunctions $\psi^*$ of the
operator \eqref{Pen.1} are  analytic (and also are analytic at
infinity, in a certain sense). However, this is out of the scope of this work. 

\vspace{0.2cm}

 \noi{\bf Remark.} The differential part in \ef{Pen.1}
can be reduced to a symmetric form in a weighted
$L^2_{\rho_\l}$-metric:
  \be
  \label{symm12}
  \tex{
  (1+z^2) D_z^2+ 2(\l+1) z D_z \equiv (1+z^2) \frac
  1{\rho_\l} D_z(\rho_\l D_z) \whereA \rho_\l =(1+z^2)^{\l+1}.
   }
   \ee
 Note that this weighted metric has an essential dependence
 on the unknown {\em a priori} eigenvalues. However, for a fixed
 $\l=\l_l$, we will use later the symmetric form \ef{symm12} in our
 branching analysis of the $p$-Laplacian problem \ef{pa}.

\vspace{0.2cm}

Currently, since our pencil approach is
nothing more than re-writing via scaling the standard
Sturm--Liouville eigenvalue problem for harmonic polynomials, it
quite natural  to deal with nothing else than them, which, thus,
should be re-built in terms of the scaling variable $z$. 

Indeed, our pencil eigenvalue problem \eqref{Pen.1} admits a reduction to a Sturm--Liouville problem whose eigenfunctions 
are harmonic polynomials.
It is easy to see that, e.g., this can be achieved by the
transformation
 \[
 \psi^*(z)=(1+z^2)^\g \varphi(z),
  \]
  with a
parameter $\g \in \re$ to be determined. Then, we find that the
operator \eqref{Pen.1} can be written as
\be
\label{stuli}
\begin{split}
 (1+z^2)^\g [ \l (\l+1)\varphi &
 +4(\l+1) \g z^2 (1+z^2)^{-1}\varphi +2(\l+1) z \varphi'+2\g \varphi  \\ &+
4z^2\g(\g-1)(1+z^2)^{-1} \varphi  +4z \g \varphi'
+(1+z^2)\varphi'']=0,
\end{split}
\ee
since
\[
(\psi^*)'(z)=2\g z (1+z^2)^{\g-1}\varphi(z)+(1+z^2)^{\g}\varphi'(z),
 \quad \mbox{and}
 \]
\[
(\psi^*)''(z)= 2\g  (1+z^2)^{\g-1}\varphi(z)+4\g (\g-1) z^2
(1+z^2)^{\g-2}\varphi(z)+ 4\g z (1+z^2)^{\g-1}\varphi'(z)+
(1+z^2)^{\g}\varphi''(z).
  \]
   To eliminate the necessary terms in
order to get a Sturm--Liouville problem, we have to cancel the
term containing $z \var'$, i.e., to require
   \[
 \tex{
  2
  \l+4\g+2\LongA
\g=-\frac{\l+1}{2}.
 }
 \]
  Now, rearranging terms for that specific $\g$ in the equation
\eqref{stuli}, so that the terms with $\varphi$ are given by
 \[
  \tex{\left(1-z^2(1+z^2)^{-1}\right)(\l+1)(\l-1)\varphi \equiv
  (1+z^2)^{-1} (\l+1)(\l-1),
  }
  \]
we arrive at a Sturm--Liouville problem of the form
\be
\label{StumLiv}
 \tex{
\mc{A} \varphi=\mu \varphi,\quad \hbox{where} \quad
\mc{A}=-(1+z^2)^2 \frac{{\mathrm d}^2}{{\mathrm d}z^2} \andA
\mu=(\l+1)(\l-1),
 }
\ee
  in the space of functions
  \[
 \tex{
  \mc{D}=L^2\big(  \re,\frac{{\mathrm d}z}{(1+z^2)^2}\big).
 }
  \]
 The
operator $\mc{A}$ is symmetric in a weighted $L^2$-space, so the
eigenvalues $\mu$ are real and by classic Sturm--Liouville
theory we state the following (see \cite{Det,Fun} for further details).

\begin{lemma}
The Sturm--Liouville problem \eqref{StumLiv} possesses a countable family 
of eigenpairs $\{\mu_n,\varphi_n\}$, such that each eigenfunction $\varphi_n$ is 
associated with the eigenvalue $\mu_n$ and the discrete family of eigenvalues satisfies
\be
\label{point}
\mu_1<\mu_2<\cdots<\mu_n \to \infty.
\ee 
Moreover, the eigenfunctions $\varphi_n$ have exactly
$n-1$ zeros in $\re$ and are the so-called $n$-th fundamental
solution of the Sturm--Liouville problem \eqref{StumLiv}. This eigenfunctions also form
an orthogonal basis in a specific weighted $L^2$-space, denoted by
$L^2_\rho$ for an appropriate weight $\rho=(1+z^2)^{-2}$.
\end{lemma}

\vspace{0.2cm}

 \noi{\bf Remark.} By classical spectral theory we also have 
that the first eigenvalue $\mu_1$ is positive and, hence thanks to  \eqref{point}
all the others. 

Also, since the weight $\rho$ is integrable, i.e.
\[\tex{\int_{\re} \frac{{\mathrm d}z}{(1+z^2)^2} <\infty,}\]
by classical spectral theory it follows that the spectrum is formed by a 
discrete family of eigenvalues as well. 
Thus, our pencil eigenvalues are associated with standard $\mu$'s via
 the quadratic algebraic equation
 \[
 \mu=(\l+1)(\l-1),
 \]
 and the correspondence of eigenfunctions
 is straightforward.

\vspace{0.2cm}

Hence, since the eigenfunctions are harmonic polynomials and, as usual in orthogonal polynomial theory, we can now state the
following property for the eigenfunctions of the adjoint pencil \eqref{Pen.1} with respect to its family of eigenvalues that will be determined below;
see \cite{Det,Fun} for details about this Sturm--Liouville Theory.

\begin{proposition}
 \label{Pr.PolP}
 The only acceptable eigenfunctions of the adjoint pencil $(\ref{Pen.1})$
are finite polynomials.
 \end{proposition}

 \vspace{0.2cm}

Although we cannot forget that once the rescaling \eqref{resc1} is performed, these eigenfunctions of the quadratic pencil operator $(\ref{Pen.1})$ 
are actually harmonic polynomials (just introducing the variables \eqref{resc1}) for which it is well known 
they are finite polynomials, it should be pointed out that this is associated with the interior elliptic
regularity. 

Indeed, the blow-up approach under the rescaling \eqref{resc1} just
specifies local structure of multiple zeros of analytic functions
at 0, and since all of them are finite (we are assuming \eqref{str1} with a finite number of cracks) we must have finite
polynomials only. 

Of course, there are other formal eigenfunctions
(we will present an example; see \eqref{001} below), but those, in the limit as $\t \to +
\iy$ in \ef{YYY.112S}, lead to non-analytic (or even discontinuous)
solutions $u(x,y)$ at 0, that are non-existent.

 \vspace{0.2cm}
 
 Moreover, the next lemma shows 
 the corresponding point spectrum of the pencil \eqref{Pen.1}. 
 
 \begin{lemma}
 The quadratic pencil operator \eqref{Pen.1} admits two families of eigenfunctions 
 $$
  \tex{\psi_{l_+}^*(z)\equiv \psi_{l,1}^*(z) \quad \hbox{and}
  \quad \psi_{l_-}^*(z)\equiv \psi_{l-1,2}^*(z), \quad \hbox{for any}\quad  l=m,m+1,\cdots}
  $$
associated with two corresponding families of eigenvalues 
 \be
 \label{ei11}
 \tex{
  \l_l^+=-l, \quad l=1,2,3,... \andA \l_l^-=-l-1, \quad
  l=0,1,2,3,... \,,
  }
  \ee
 \end{lemma}

\noi{\em Proof.}
 In order to find the corresponding point spectrum of the pencil we look for $l$th-order polynomial eigenfunctions of the form
 \be
 \label{eig00}
  \tex{
\psi_l^*(z)=z^l+a_{l-2}z^{l-2}+ a_{l-4}z^{l-4}+... =
\sum\limits_{k=l,l-2,...,0} a_k z^k \quad (a_l=1)
  }
   \ee
   that we already know they are harmonic polynomials. Substituting \eqref{eig00} into (\ref{Pen.1})
   and evaluating the higher order terms
 yields the following quadratic equation for eigenvalues:
  \be
  \label{eig11}
 O(z^l): \quad  \l_l^2 + (2l+1) \l_l + l(l + 1)=0.
   \ee
Solving this characteristic equation yields the two families of real
negative eigenvalues under the expression \eqref{ei11}
  associated with two families of eigenfunctions denoted by
 \be
 \label{harpol}
  \tex{\psi_{l_+}^*(z)\equiv \psi_{l,1}^*(z) \quad \hbox{and}
  \quad \psi_{l_-}^*(z)\equiv \psi_{l-1,2}^*(z), \quad \hbox{for any}\quad  l=m,m+1,\cdots}
  \ee
for convenience.  \quad $\qed$

 \vspace{0.2cm}
 
 The next result calculates those  (re-structured harmonic)
  polynomials \eqref{eig00}, \eqref{harpol} as the corresponding eigenfunctions of the pencil. 

 \begin{theorem}
  \label{Th.1}
 The quadratic pencil \ef{Pen.1} has two (admissible) discrete spectra
 $\ef{ei11}$ of real negative eigenvalues with the finite polynomial
 eigenfunctions given by \ef{eig00}, where the expansion
 coefficients satisfy a
 finite Kummer-type recursion corresponding to the operator in
$(\ref{Pen.1})$:
  \be
  \label{co1}
   \left\{
    \begin{array}{ccc} \tex{
    a_{k+2}= -
    \frac{k(k-1) +2(\l_l^\pm+1)k +\l_l^\pm(\l_l^\pm+1)
   }{(k+2)(k+1)} \,  a_{k},} & \hbox{for any} & k=l,l-2,...,2,
         \\
        \tex{ a_1= -\frac{6}{ 2(\l_l^\pm+1)  +\l_l^\pm(\l_l^\pm+1)}a_3} & \hbox{and} & \tex{a_0 =- \frac{2}{\l_l^\pm(\l_l^\pm+1)} a_2.}
        \end{array}\right.
    \ee
  \end{theorem}

\noi{\em Proof.}
    It is clear by \eqref{ei11} that the quadratic pencil \ef{Pen.1} has two discrete spectra
    of real negative eigenvalues with two families of finite polynomial
 eigenfunctions\footnote{Note that, within this pencil ideology,
 the eigenfunctions are ordered in an unusual manner, unlike
 the standard harmonic polynomials.}
$$
\{\psi_{l_+}^*(z)\},\quad \{\psi_{l_-}^*(z)\},\quad \hbox{such that}\quad \psi_{l_+}^*(z)
\equiv \psi_{l,1}^*(z)\quad \hbox{and}\quad \psi_{l_-}^*(z)\equiv \psi_{l-1,2}^*(z),
$$
given by \ef{eig00} and corresponding associated with the two families of eigenvalues $\l_l^+$ and $\l_l^{-}$ .
Substituting $\psi_l^*=\sum\limits_{k\geq 0}^l a_k z^k$, for any $l\geq 0$, into (\ref{Pen.1})
    we find that, for any $\l$,
$$
 \tex{
 (1+z^2) \sum\limits_{k\geq 2}^l k(k-1) a_k z^{k-2} + 2(\l+1) \sum\limits_{k\geq 1}^l k a_k
z^{k}+\l(\l+1) \sum\limits_{k\geq 0}^l a_k z^{k}=0,
 }
 $$
and hence,
\be
\label{coefla}
\begin{split}
\sum\limits_{k\geq 2}^l & \left[(k+2)(k+1) a_{k+2}+ k(k-1) a_k+2(\l+1)ka_k +\l(\l+1)a_k\right] z^k \\ & +
\left[ 6a_3 + [2(\l+1)  +\l(\l+1)]a_1\right] z+ 2a_2+\l(\l+1) a_0=0.
\end{split}
\ee
Therefore, evaluating the coefficients we find that
$$\left\{\begin{array}{l}
    (k+2)(k+1) a_{k+2}+ k(k-1) a_k+2(\l+1)ka_k +\l(\l+1)a_k=0,\quad k=l,l-2,...,2,\\
    6a_3 + [2(\l+1)  +\l(\l+1)]a_1=0,\\
    2a_2+\l(\l+1) a_0=0,\end{array}\right.
    $$
    and we arrive at  \eqref{co1}, completing the proof. \quad $\qed$

 \vspace{0.2cm}

Note that even when
discrete spectra coincide excluding the first eigenvalue
$\l_l^-$, and, more precisely,
$$
\l_l^-= \l_l^+-1= \l_{l-1}^+ \quad l=1,2,3,...\,,
$$
we still have two different families of eigenfunctions. For future
convenience and applications for the crack problem for $m=1,2,3$,
and 4 (with $m=l$), we present first four eigenvalue-eigenfunction
pairs of both families of eigenfunctions for the pencil
\ef{Pen.1}, which now are ordered with respect to $\l=-l$,
$l=0,1,2,...$:
\be
\label{pol111}
 \begin{split}
\l_0=0, &  \quad \mbox{with} \,\,\, \psi^*_0(z)\equiv
\psi_{0,1}^*=1 \,\,\,(\ne 0); \\ \l_1=-1, & \quad \mbox{with}
\,\,\, \psi^*_{1,1}(z)=z, \quad \psi^*_{0,2}(z)=1 \,\,\,(\ne 0);\\
 \l_2=-2, & \quad \mbox{with} \,\,\,\psi^*_{2,1}(z)=z^2-1,
 \,\,\, \psi^*_{1,2}(z)=z; \\
  \l_3=-3, & \quad \mbox{with} \,\,\,
 \psi^*_{3,1}(z)=z^3 - 3z,\,\,\,
\psi^*_{2,2}(z)=3z^2-1; \\
 \l_4=-4, & \quad \mbox{with} \,\,\,
\psi^*_{4,1}(z)=z^4 - 6z^2+1,\,\,\,
 \psi_{3,2}^*(z)=z^3 -z;\,\,\, \mbox{etc.}
 \end{split}
 \ee
 
 \vspace{0.2cm}

\noi{\bf Remark: about transverality.}
These (harmonic) polynomials satisfy the Sturmian property (important for applications) in the sense that
 each polynomial $\psi_{m}^*(z)$ has precisely $m$ {\em transversal}
 zeros. 
 For Hermite polynomials, this result
  was proved by Sturm already in 1836
  \cite{St}; see further historical comments in
  \cite[Ch.~1]{GalGeom}.

 \vspace{0.2cm}

\noi{\bf Remark:
about analyticity.} Obviously, we exclude,
in the first line of \ef{pol111}, the first eigenfunction
$\psi^*_0(z) \equiv \psi^*_{0,1}(z) \equiv 1$, since it does not
vanish and has nothing to do with a multiple zero formation.
However, for $\l=0$ in \ef{Pen.1}, there exists another obvious
bounded analytic solution having a single zero:
 \be
 \label{001}
 (1+z^2) (\psi^*)''+ 2z (\psi^*)'=0 \LongA \tilde\psi^*(z)= \tan^{-1}z \to \pm \pi/2 \asA z \to \pm \iy.
 \ee
This $\tilde \psi^*(z)$ belongs to any suitable $L^2_\rho$-space
(of polynomials). However, it becomes irrelevant due to another
regularity reason: passing to the limit in the corresponding
expansion of $u(x,y) \equiv w(y,\t)$ \ef{YYY.112S} as $\t \to
+\iy$ ($y \to -0$) yields the discontinuous limit ${\rm sign}\,x$,
i.e., an impossible trace at $y=0$ of any analytic solutions of
the Laplace equation.

\subsection{Nonexistence result for the crack problem \eqref{Lap1}
}

Next we ascertain how the family of admissible cracks should lead to the existence of solutions for the crack problem \eqref{Lap1}.

It is well known that sufficiently ``ordinary" polynomials are
always complete in any reasonable weighted $L^2$ space, to say
nothing about the harmonic ones; see \cite[p.~431]{KolF}.
Moreover, since our  polynomials are not that different from
harmonic (or Hermite) ones, this implies the completeness
in such spaces. 

So that, sufficiently regular solutions of
\ef{www.2}
should admit the corresponding
eigenfunction expansions over the polynomial family pair $\Phi^*=
\{\psi_{l,1}^*,\psi_{l-1,2}^*\}$ in the following sense. 
\begin{lemma}
Bearing in mind two discrete
spectra \ef{ei11} and \eqref{YYY.112S}, the general expansion for the solutions of \eqref{www.2}
has the form
 \be
 \label{exp11}
  \tex{
  w(z,\t)= \sum_{(k \ge l)} \eee^{-k \t} [c_k\psi_{k,1}^*(z)+ d_k \psi_{k-1,2}^*(z)],
  }
  \ee
 where two collections of expansion coefficients $\{c_k\}$ and
 $\{d_k\}$,
  depending on boundary
 data on $\O\setminus \G$, are presented. 
 \end{lemma}
 
 To develop an
 ``orthonormal theory" of our polynomials, we should specify the
 expansion coefficients in \ef{exp11}, for a given solution
 $u(x,y)$ (though specifying all the coefficients declare the
 whole family of $u$ with such cracks at 0). Then, one just can transform 
 the standard expansion for the harmonic solutions (orthogonal harmonic polynomials) of the Laplace problem 
 and obtain \ef{exp11} by introducing the scaling
 blow-up variables \ef{resc1}.

 Furthermore, the linear combination \eqref{exp11} arises naturally
from the spectral theory of the operator (in this case the Laplacian, later on the bi-Laplacian).
Indeed, for the Laplacian $u$ is harmonic in $B_1\setminus \G$ and can be decomposed
by homogeneous harmonic functions, here denoted by $\psi_{k,1}^*$ and $\psi_{k-1,2}^*$. Even facing
a difficult regularity problem in $\O\setminus \G$ (at the singularity boundary point) we are in the context analysed in
\cite{AL}, so that
 $$
 \tex{w(z,\t)=  \eee^{-k \t} \psi^*(z),}
 $$
 for an orthonormal basis $\{\psi_{k,1}^*,\psi_{k-1,2}^*\}$ of Hermite-type polynomials eigenfunctions. Hence, we find that
 our solutions are decompositions of the form \eqref{exp11}.

 Moreover, in view of sufficient regularity of ``elliptic
orbits" (via interior elliptic regularity), such
expansion is to converge not only in the mean (in $L^2_\rho$,
with an exponentially decaying weight at infinity), but also
uniformly on compact subsets. This allows us now to prove our
result on nonexistence for the crack problem.

\begin{theorem}
\label{Th.2} Let the cracks $\G_1$,...,$\G_m$ in \ef{curv1} be
asymptotically given by $m$ different straight lines \ef{str1}.
Then, the following hold:

{\rm (i)} If all
$\{\a_k\}$ do not coincide with  all $m$ subsequent zeros of any
non-trivial linear combination
 \be
 \label{dmex}
c_l \psi^*_{l,1}(z)+d_l \psi^*_{l-1,2}(z), \quad \mbox{with} \quad
c_l^2+d_l^2 \ne 0 \whereA z= x/(-y),
 \ee
 of two families of (re-written harmonic)
polynomials $\psi_{l_+}^*(z)\equiv  \psi^*_{l,1}(z)$ and $\psi_{l_-}^*(z)\equiv  \psi^*_{l-1,2}(z)$ defined by
\ef{eig00}, \ef{co1} for any $l=m,m+1,...$ and arbitrary constants
$c_l, \, d_l \in \re$, then the multiple crack problem \ef{Lap1}
cannot have a solution for any boundary Dirichlet data $f$ on
$\O$.

{\rm (ii)} If, for some $l$, the distribution of zeros in {\rm
(i)} holds and a solution $u(x,y)$ exists,
 then
 \be
 \label{as11ex}
  \tex{
  |u(x,y)| = O(|x,y|^{l})    \asA (x,y) \to (0,0).
  }
  \ee

 \end{theorem}

 \noi{\em Proof of Theorem \ref{Th.2}.} Condition \ef{str1} implies that the elliptic
 ``evolution" problem while approaching the origin actually occurs on
 compact, arbitrarily large subsets for $ x/(-y) \equiv z$.
 Since we have converted the singularity point at $(0,0)$ into an
 asymptotic point when $\tau \to \infty$.

 Therefore, \ef{exp11} gives all possible types of such a decay.
 Hence, choosing the first non-zero expansion coefficients $c_l$, $d_l$ in
 \ef{exp11}, that satisfies $c_l^2+d_l^2 \ne 0$, we obtain a sharp asymptotic behaviour of this solution
  \be
  \label{st.2}
  w_l(y,\t)=  \eee^{-l \t}[c_l \psi_{l,1}^*(z)+d_{l} \psi_{l-1,2}^*] + O(\eee^{-(l+1) \t})
  \asA \t \to +\iy.
  \ee
  Obviously,  
then the straight-line cracks \ef{str1} correspond to zeros of
the linear combination 
$$\tex{c_l \psi_{l,1}^*(z)+d_{l}
\psi_{l-1,2}^*(z)},$$ 
and the full result is straightforward since by the blow-up scaling if all the $\a_k$ do not 
coincide with zeros of the previous linear combination \eqref{dmex} (harmonic polynomials) the crack problem 
does not have a solution, since
$$ z=\frac{x}{-y}=\a_k (1+o(1)), \,\, y \to 0.$$
 Otherwise, if there is some $l$ for which all the $\a_k$
coincide with zeros of \eqref{dmex} we find that the crack problem \eqref{Lap1} possesses a solution and \eqref{as11ex} is satisfied.
  The proof is complete. \quad
$\qed$

\vspace{0.2cm}

\noi{\bf Remark.}
Of course, one can ``improve" such nonexistence results. For
instance, if cracks have an asymptotically small ``violation" of
their straight line forms near the origin, which do not correspond to
the exponential perturbation in \ef{st.2} (if $c_{l+1}$ and
$d_{l+1}$ do not vanish simultaneously; otherwise take the next
non-zero term), then the crack problem is non-solvable.

\ssk

Overall, we can state the following most general conclusion. 

\begin{corollary}
For almost every straight-line crack \ef{curv1}, the crack problem
\ef{Lap1} cannot have a solution for any Dirichlet data $f$,
provided that the crack behaviour at the origin is not consistent
with all the eigenfunction expansions \ef{exp11} via the above
(harmonic) polynomials.
\end{corollary}

\ssk

Finally, concerning the admissible boundary data for such
$l$-cracks at the origin, these are described by all the
expansions \ef{exp11} with arbitrary expansion coefficients
excluding the first ones $c_l$, $d_l$, which are fixed by the
multiple crack configuration (up to a common non-zero multiplier)
and satisfying $c_l^2+d_l^2 \ne 0$.



\section{$p$-Laplace equation: nonlinear eigenfunction and branching at $n=0$}
\label{plaplace}

\subsection{Rescaled ``evolution" equation}

We now consider the $p$-Laplace equation \ef{pa} in the equivalent and
more convenient  form, for our analysis,
 \be
 \label{A1}
 \D_p u \equiv [(u_x)^2+(u_y)^2]^{\frac n2} \D u+ n [(u_x)^2+(u_y)^2]^{\frac {n-2}2}
 [(u_x)^2 u_{xx}+ 2 u_x u_y u_{xy}+(u_y)^2 u_{yy}]=0,
  \ee
  with $n=p-2$.
Using the same blow-up variables \ef{resc1} as those ones used for the Laplace problem \eqref{Lap1} in the 
previous section yields
 \be
 \label{A2}
  \begin{aligned}
  &\D_p u \equiv  [(w_z)^2+(w_\t+z w_z)^2]^{\frac n2} \D_{(z,\t)}w + n [(w_z)^2+(w_\t+zw_z)^2]^{\frac {n-2}2}
  \big[(w_z)^2 w_{zz}\\
  & +  \,2w_z (w_\t+z w_z)(w_{z \t}+ z w_{zz})+ (w_\t+z
  w_z)^2(w_{\t\t} + w_\t +2 z w_{z \t}+ z^2 w_{zz} + 2 z w_z)\big]=0.
   \end{aligned}
   \ee
Thus, choosing again $\t$ as the evolution variable  and using the
Laplace operator in \ef{www.20}, instead of \ef{www.2}, we arrive
at a quasilinear equation
 \be
 \label{A3}
  \begin{aligned}
   &
   \tex{(w_{\t\t}+w_\t+2 z w_{z \t}+z^2 w_{zz}+ 2z w_z)\big[
  1+ \frac {n(w_\t+zw_z)^2}{(w_z)^2+(w_\t+zw_z)^2} \big]}
   \\
  &
  \tex{ = \, - w_{zz} - \frac{n[(w_z)^2 w_{zz}+2w_z(w_\t+z w_z)(w_{z \t}+z
  w_{zz})]}{(w_z)^2+(w_\t+zw_z)^2}.
  }
  \end{aligned}
 \ee
In particular, for $n=0$, \ef{A3} formally coincides with \ef{www.2}. We actually choose the more convenient form \eqref{A1} to get such a 
relation.

\subsection{Non-linear eigenvalue problem}

 The non-linear PDE
\ef{A3} remains homogeneous of degree 1, i.e., if $w$ is a
solution, then $C w$ is also  one, for any constant $C \in \re$.
Hence, it admits separation of variables as in
 \ef{YYY.112S}:
 \be
 \label{A40}
 w(z,\t)=\eee^{\Lambda \t} \Psi^*(z),
  \ee
  that leads to a {\em non-linear
 eigenvalue problem}
 \be
 \label{A4}
 \tex{
  [\Lambda(\Lambda+1) \Psi+ (2 \Lambda+1) z \Psi'+z^2\Psi'']
  [1+n \Phi_1(\Psi,\Lambda)]=-\Psi''- n \Phi_2(\Psi,\Lambda), \quad
  \mbox{where}
  }
  \ee
  \be
  \label{A5}
   \tex{
   \Phi_1(\Psi,\Lambda)= \frac{(\Lambda \Psi+z
   \Psi')^2}{(\Psi')^2+(\Lambda
   \Psi+z \Psi')^2}
   \andA \Phi_2(\Psi,\Lambda)= \frac{(\Psi')^2\Psi''+2\Psi'(\Lambda \Psi+z
   \Psi')(\Lambda \Psi'+z \Psi'')} {(\Psi')^2+(\Lambda
   \Psi+z \Psi')^2},
    }
    \ee
    for which we intend to  find all {\em real}\footnote{Cf. real ``linear" eigenvalues \ef{ei11} of the pencil
 \ef{Pen.1}.} eigenvalues
 $\Lambda<0$ so that there exists an ``admissible" (see below) nonlinear
 eigenfunction $\Psi$ (for the moment, we omit the superscript $*$).

   Furthermore,  to complete this {\em nonlinear eigenvalue problem} for \ef{A4},
    one needs proper singular ``boundary conditions at infinity",
    and this is not that easy in such a nonlinear setting.
However, it turns out that we can use here quite similar
conditions as in the above linear case $n=0$. 

\begin{lemma}
Non-linear eigenfunctions $\Psi$ possess a polynomial growth at infinity of the form 
\be
      \label{A6}
 \Psi^*_l(z)=z^l(1+o(1)) \asA \to \iy,
  \ee
  if they are uniformly non-degenerate at infinity.
\end{lemma}

\noi{\em Proof}. Recall that, for
    $n=0$ (see \ef{eig00}), we allow only a polynomial growth of
    nonlinear eigenfunctions, i.e., a polynomial growth of the type \eqref{A6} 
  as a linear combination of the two families
 $$\{\psi_{l_+}^*(z)\},\quad \{\psi_{l_-}^*(z)\},\quad
\hbox{such that}\quad \psi_{l_+}^*(z) \equiv \psi_{l,1}^*(z)\quad
\hbox{and}\quad \psi_{l_-}^*(z)\equiv \psi_{l-1,2}^*(z).$$
 At $n=0$ (for the Laplace problem \eqref{Lap1}) this was connected with the analyticity of solutions, so that each zero at
 the origin should be of a finite order $l \ge 1$, where harmonic
 polynomials locally and the extension \ef{A6} appeared from. Indeed finite harmonic polynomials are responsible
 for all types of local behaviour of multiple zeros around 0. Remember that we have converted the singularity point $(0,0)$ 
 into an asymptotic convergence when $\tau\to \infty$. 

 Thanks to elliptic interior regularity, solutions of \ef{A1} and, hence \ef{A4}, are principally
 non-analytic and have finite regularity, but {\em at points of
 degeneracy only}, where 
 $$\Psi'(z)=0.$$ This is easily spotted from the expression of the equation \eqref{A4} and 
 \eqref{A5}. 
 
 Beyond non-degeneracy sets, solutions of the
 ODE problems like \ef{A4} {\em are analytic} by classic ODE
 theory. Basically due to the fact that the singular terms \eqref{A5} are not singular any more. 
 
 Therefore, the polynomial growth \ef{A6} remains in
 charge for $n>0$ as well, provided that the corresponding nonlinear
 eigenfunctions
 $$\Psi(z)\quad \hbox{are uniformly non-degenerate at
 infinity, i.e., $|\Psi'(z)| \ge \d_0>0$ for all $|z| \gg 1$}.$$  $\qed$
 
 \vspace{0.2cm}

 \noi{\bf Remark.} Note that with such a ``linear" condition at infinity \ef{A6}
can be also associated with the above mentioned 1-homogeneity of
nonlinear operators in \ef{A3}.  In other words, due
to those ``linear" properties of these operators, the nonlinear
eigenvalue problem inherits the ``linear" conditions (normalization)
\ef{A6}. 

Furthermore, one can see that exactly this 1-homogeneity of
operators  allows us to obtain a ``nonlinear characteristic equation".
   Thus,
  substituting \ef{A6} into \ef{A4} yields the following
  ``nonlinear characteristic polynomial equation" for eigenvalues $\{\Lambda_l\}_{l
  \ge 1}$:
 \be
   \label{A7}
     \tex{
 [\Lambda_l^2 +(2l+1)\Lambda_l + l(l+1)]\big[1+
 n \, \frac {(\Lambda_l+l)^2}{l^2+(\Lambda_l+l)^2} \big]}
  \tex{ + \,n \, \frac{l^3(l-1)+2l(\Lambda_l +l)(\Lambda_l l +l(l-1))}{l^2+(\Lambda_l +l)^2}=0.
  }
 \ee
For $n=0$, \ef{A7} leads to the quadratic equation
 \ef{eig11}.
One can see that \ef{A7} reduces to the following quartic
``characteristic" equation:
 \be
 \label{A71}
  \begin{matrix}
  \Phi_l(\Lambda;n) \equiv (1+n)\Lambda^4 + a_3 \Lambda^3+ a_2 \Lambda^2 + a_1 \Lambda +
  a_0, \,\,\, \mbox{where} \ssk \\
 a_3=(1+n)(4l+1),\,\,\,
  a_2=l[5n+3+l(6n+7)],\ssk
 \\
  a_1=l^2[3n+4+6l(1+n)], \,\,\,
  a_0=l^3[2(1-n)+2l(2n+1)].
  \end{matrix}
  \ee
 Once the nonlinear algebraic eigenvalue equation has been
 solved and real negative admissible nonlinear spectrum $\{\Lambda_l\}$
 has been obtained (this can be done analytically for small $n>0$ with
 numerical extensions; see below), one arrives at the
 corresponding nonlinear eigenfunction problem: given a real value
 $\Lambda_l<0$, to find a nontrivial solution $\Psi=\Psi^*_l(z) \not
 \equiv 0$ of the nonlinear ODE
   \be
   \label{A9}
    \left\{
 \begin{aligned}
  &
     \tex{
 [\Lambda_l(\Lambda_l+1)\Psi + 2 (\Lambda_l+1) z\Psi'+z^2 \Psi'']\big[1+
  n\,\frac {(\Lambda_l \Psi+z \Psi')^2}{(\Psi')^2+(\Lambda_l \Psi+z \Psi')^2} \big]}
   \\
  &
  \tex{ = \, - \Psi'' - n\,  \frac{(\Psi')^2 \Psi''+2\Psi'(\Lambda_l \Psi+z \Psi')(\Lambda_l \Psi'+z
  \Psi'')}{(\Psi')^2+(\Lambda_l\Psi+z \Psi')^2}, \quad z \in \re,
  } \\
 &
 \mbox{$\Psi=\Psi^*_l(z)$ has the polynomial growth (\ref{A6}) at infinity}.
  \end{aligned}
  \right.
 \ee
 \vspace{0.2cm}

 \noindent\underline{Non-linear eigenfunctions $\Psi$ with transversal zeros}.
Equation \eqref{A9} looks rather frightening, but, in fact, for any
solution $\Psi(z)$ not having {\em non-transversal zeros}, where,
for some $z_0>0$,
 \be
 \label{psi99}
\Psi(z_0)=\Psi'(z_0)=0\quad \hbox{(condition for non-transversal
zeros)}
 \ee
  (actually, this never happens; see below). Both
denominators in \ef{A9} do not vanish, so the solution is analytic
(as in the linear case $n=0$). So that, in the limit as $\t \to
+\iy$, it  must be represented by a zero of $u(x,y)$ at the origin
with finite multiplicity. Thanks to the polynomial condition \eqref{A6}. 

Concerning the non-transversal zeros
\ef{psi99}, the transversality condition follows from an easy local analysis of the ODE
\ef{A9} near $z=z_0$. 

\begin{lemma}
Non-linear eigenfunctions $\Psi$ satisfy a transversality condition in all their zeros.
\end{lemma}

\noi{\em Proof}. Thanks to \ef{psi99}, we can assume that
$$
 \tex{\Psi(z) =O(\Psi'(z))\quad \hbox{as}\quad z \to z_0,}
   $$
   (excluding an oscillatory
behaviour for such a second-order ODE \eqref{A9}), so that both ``singular"
terms in \ef{A5} satisfy,
 \be
 \label{psi98}
  \tex{
 \Phi_1 \sim \frac {z_0^2}{1+z_0^2} \,\,\,(\mbox{bounded}) \andA
 \Phi_2 \sim \frac{z_0^2}{1+z_0^2}\, \Psi''
 \,\,\,(\mbox{uniform}) \quad \hbox{as}\quad z \to z_0.
  }
  \ee
  Just taking into consideration the polynomial behaviour \eqref{A6} and substituting them into the expressions of the singular terms \eqref{A5}. 
  
  Therefore, substituting \eqref{A6} into the equation \eqref{A9} and due to the local behaviour of the singular terms \eqref{psi98} around the zeros $z_0$ 
  we arrive at the fact that the ODE \ef{A9}, with conditions \ef{psi99}, gives
$$
 \hbox{the unique solution $\Psi=0$ for all $z$, if there is
non-transversal zeros}.
 $$
 Thus, \ef{A9} does not admit solutions
with non-transversal zeros as in \ef{psi99}. $\qed$
 
 \vspace{0.2cm}

Note that first two pairs  $\{\Psi_0^*, \Lambda_0\}$ and
$\{\Psi_1^*, \Lambda_1\}$ are easy and these  are the same as in
the linear case, see \ef{pol111}: for all $n \ge 0$,
 \be
 \label{A91}
 \Lambda_0 (n)=0, \quad \Psi_0^*(z) \equiv 1 \ne 0 \andA
 \Lambda_1(n)=-1, \quad \Psi_1^*(z)=z.
  \ee
  
   \vspace{0.2cm}

\noindent{\bf Remark.} As we have mentioned
above, for the Hermite polynomials/eigenfunctions of the
adjoint Hermite operator 
$$\BB^*= D_z^2- \frac 12 \, z D_z,\quad \hbox{in}\quad L^2_\rho,\quad \hbox{with}\quad \rho=\eee^{-z^2/4},$$ 
the transversality of all their
zeros
  was proved by Sturm already in 1836
  \cite{St}; see comments in
  \cite[Ch.~1]{GalGeom}.  Curiously, we see now that this kind of a Sturmian
  zero transversality property remains valid for the quasilinear
  problem \ef{A9}.
As a consequence, we conclude that its solutions are analytic
functions for  all $z$.

\vspace{0.2cm}

  \noi{\bf Remark: a discontinuous limit.} As for $n=0$, let us
  show how a wrong choice of an eigenfunction can contradict the
  interior regularity for the $p$-Laplacian. Setting $\Lambda=0$
  in \ef{A9} yields the following ODE, which can be integrated once:
   \be
   \label{A99}
    \tex{
   (1+n)(1+z^2)^2 \Psi''+ 2z[1+(1+n)z^2] \Psi'=0
   \LongA
   \Psi'(z)= \frac 1{1+z^2} \exp\{- \frac n{1+n} \frac 1{1+z^2}\}.
   }
   \ee
   A further implicit integration yields an ``eigenfunction"
   satisfying 
   $$\tilde \Psi(z) \to \pm 1,\quad \hbox{as}\quad  z \to \pm \iy,$$ 
   with a
   single zero at $z=0$. Then, passing to the limit as $y \to 0^{-}$ in
   \ef{A40} gives a discontinued limit ${\rm sign} \, x$, not
   suitable for the $p$-Laplacian via interior regularity.

\subsection{Existence of eigenfunctions}

 The above analysis
implies that, once $\{\Lambda_l(n)\}$ are known, existence of
eigenfunctions are straightforward, since the ODE \ef{A9} does not
admit any blow-up, so local solutions are globally extensible  for
all $z>0$. 

Therefore, with correct eigenvalues (putting proper
symmetry, for even $l$, or anti-symmetry, for odd $l$, conditions
at $z=0$), one obtains global solutions, which, inevitably,
satisfy the proper polynomial behaviour at infinity \ef{A6}. The
latter is right, since no other behaviour at infinity is available
(we omit certain technicalities establishing such local properties
of the ODE at $z=\iy$). 

In fact, there appears to be no boundary value
problem in this setting: e.g., take an even $l \ge 2$
 and, in the Cauchy problem for \ef{A9} for $z>0$, put the initial  conditions
  \be
  \label{bv43}
  \Psi(0)=C \ne 0 \andA \Psi'(0)=0.
  \ee
  This problem contains no parameters, since $C$ can be scaled out
  (and reduced to 1) by the ``linear" 1-homogeneity property of
  the operators in \ef{A9}. Hence, due to the election of conditions \ef{bv43} uniquely
  determines the eigenfunction $\Psi_l(z)$, which, by the choice
  of $\Lambda_l(n)$, satisfies the necessary polynomial behaviour
  at infinity.

\ssk

Therefore, the main difficulty of this nonlinear eigenvalue
problem lies in the study of the algebraic characteristic
equation \ef{A71}.
  However, we begin with a simpler existence analysis, simultaneously, of both such nonlinear  eigenvalues and
 eigenfunctions by using  a branching
 approach at $n=0$, where, again and again, harmonic polynomials naturally occur.


 \section{Branching at $n=0$ of nonlinear eigenfunctions
 from harmonic polynomials}
 \label{SBr}

 We now show, using a branching approach at $n=0$ (following a similar philosophy to \cite{TFE4PV}), how a
 polynomial behaviour \ef{A6} of nonlinear eigenfunctions
 is connected with that for analytic harmonic polynomials. Namely, looking at the equation \ef{A9} as
a perturbed linear ODE \ef{Pen.1} for small $n>0$, with
 perturbations of order $O(n)$. We then look for solutions of
 \ef{A9}.

\begin{theorem}
\label{mainpt}
Consider  \eqref{A9} as a perturbed equation of the linear
differential equation \eqref{Pen.1}, for small $n>0$, with
perturbations of order $O(n)$. Then the solutions of the equation
\eqref{A9} have the standard form
  \be
  \label{A10}
   \left\{
   \begin{aligned}
   &
   \Psi^*_l(z)= \psi_l^*(z) + n\, \var_l(z) + o(n),\\
   &
   \Lambda_ l=\l_l+n \, \mu_l+o(n),
    \end{aligned}
    \right.
    \ee
    where $\{\psi_l^*,\l_l\}$ are linear pairs \ef{ei11},
    $\{\var_l(z)\}$ are unknown functions, and $\{\mu_l\}$ are constants, both  to
    be determined via branching equations.
        \end{theorem}
\noi{\em Proof}. First we define the operator $L_\l$ which will play a crucial role in applying the Lyapunov--Schmidt reduction. Indeed, 
by $L_\l$ we denote the first linear operator in
\ef{A9}. Thus, due to the second expansion in \eqref{A10} it follows that
 \be
 \label{A11}
 L_\l=\l(\l+1){\rm Id}+2(\l+1)z D_z+z^2 D^2_z
 \LongA
 L_{\Lambda_l}=L_{\l_l}+n \mu_l(2 \l_l+1+zD_z)+o(n).
 \ee
 Note that, thanks to \ef{Pen.1},
  \be
  \label{A12}
  \BB^*_{\l_l}=L_{\l_l}+ D^2_z.
   \ee
Applying the expansions \ef{A10} to the problem \eqref{A4} yields
\begin{align*}
 [(\l_l & +n \, \mu_l) ((\l_l+n \, \mu_l)+1)  (\psi_l^*(z) + n\, \var_l(z)) + 2 ((\l_l+n \, \mu_l)+1) z(\psi_l^*(z) + n\, \var_l(z))'
 +z^2 (\psi_l^*(z) + n\, \var_l(z))''] \\ & [1+n \Phi_1(\psi_l^*(z) + n\, \var_l(z),\l_l+n \, \mu_l)]
  =-(\psi_l^*(z) + n\, \var_l(z))''- n \Phi_2(\psi_l^*(z) + n\, \var_l(z),\l_l+n \, \mu_l) +o(n).
\end{align*}
Now, passing to the limit as $n\to 0^+$ we find that
\be
\label{A113}
 \quad \BB^*_{\l_l} \psi^*_l=0.
 \ee
  Obviously, equation in \ef{A113} coincides with the
 linear one \ef{Pen.1}, so that there exists
  \be
  \label{A14}
 \mbox{\rm ker} \, \BB_{\l_l}^*={\rm Span}\, \{\psi_l^*\}.
  \ee
Subsequently, dividing the rest of terms by $n$ and passing to the limit as $n\to 0^+$ yields
\be
\label{A13}
 \quad \BB_{\l_l}^* \var_l=h(\mu_l,\psi_l^*, \l_l)
   \equiv
 -[\Phi_2(\psi_l^*,\l_l)+\mu_l((2\l_l+1)\psi_l^*+z
 (\psi_l^{*})')+\Phi_1(\psi_l^*,\l_l)\, L_{\l_l} \psi_l^*].
\ee
  Using the symmetric form \ef{symm12} of $\BB_{\l_l}^*$ and re-writing equation
  \ef{A13} in the form
 \be
 \label{A15}
  \tex{
  \frac 1{\rho_{\l_l}}\, (\rho_{\l_l} \var_l')' +
  \frac{\l_l(\l_l+1)}{1+z^2}\, \var_l= \frac{h(\mu_l,\psi_l^*,\l_l)}{1+z^2},
   }
   \ee
   we then
  arrive at the following orthogonality condition for the unknown
  parameters:
   \be
   \label{A16}
    \tex{
   \mu_l: \quad  \int\limits_{\re} \rho_{\l_l}(z) \, \frac{h(\mu_l,\psi_l^*(z),\l_l)\,
    \psi_l^*(z)}{1+z^2}\,\, {\mathrm d} z=0
    \quad(\rho_{\l_l}(z)=(1+z^2)^{\l_l+1}).
    }
    \ee
 Since, by \ef{ei11}, 
$$\tex{\l_l=-l \le -1\quad\hbox{for}\quad l=1,2,3,...,}$$ 
integrals
in \ef{A16} make usual sense and are convergent.

Note that \eqref{A16} is obtained simply applying the Fredholm alternative after multiplying the equality \eqref{A15} by $\psi_l^*$ and integrating. 
    Indeed, applying FredholmÕs theory \cite{Deim} to \eqref{A15} yields that there exists a function $\var_l$ 
    which solves \eqref{A15} if and only if the right hand side is orthogonal to $ \mbox{\rm ker} \, \BB_{\l_l}^*$, 
    i.e., to the eigenfunction $\psi_l^*$ of the operator $\BB_{\l_l}^*$. Hence, the linear algebraic equation
\ef{A16} allows us to get unique values of $\mu_l$ \eqref{A16}.

Consequently, under condition \ef{A16}, functions $\var_l(z)$ are uniquely
determined from the second equation in \ef{A13}, since
$\Psi^*_l(z)$ defined by \eqref{A10} emanate uniquely at $n=0$ from
$\psi_l^*(z)$. $\qed$

\vspace{0.2cm}

We conclude that the nonlinear eigenfunctions
$\{\Psi_l^*(z)\}$, constructed in such a way, at least for small
$n>0$,  via their nodal sets,
allow us to specify admissible crack distributions in the
$p$-Laplacian problem \ef{pa}, \ef{curv1}, \ef{str1}.
A global continuation of such $n$-branches of
eigenfunctions requires more difficult mathematics and even
numerics; see below.

However, for a complete classification of such a crack
configuration, an {\em evolution completeness} of this set
 $\{\Psi_l^*(z)\}_{l \ge 1}$ is necessary, just meaning that it
 contains all possible limits as $x \to 0$ of solutions of \ef{pa} that
 vanish at the origin.   This is also a very difficult open problem,
 which, actually, was solved just for a couple of  much easier
 nonlinear evolution
 problems; see notions, results, and references in \cite{GalC}.

\vspace{0.2cm}

  \noi{\bf Remark.} Note that a kind of completeness of $\{\Psi_l^*(z)\}$ can be
expected, via $n$-branching, for small $n>0$, in view of
completeness/closure of harmonic polynomials for $n=0$ (though,
this should be proved: we do not know any relation between
completeness in the standard linear and nonlinear cases). However,
for larger $n>0$, we face another difficulty that we will deal with
below.


In addition, let us mention that, in this nonlinear case, in view
of the absence of any eigenfunctions expansion representation, a full
classification of boundary data that leads to such multiple zeros
of solutions at the origin, becomes intractable. In general, such
problems are solved by a matching (or extension) of asymptotic
expansions at 0 and close to $\partial \O$, but we do not believe
that, even using the 1-homogeneity of operators, any matching like
that, being always {\em asymptotic} in the nature, can specify
such sufficiently sharp data. 

Since a proper analytic approach for this global continuation analysis is still unknown 
we show here some results via numerical evidences.

\section{Towards global continuation of eigenvalue $n$-branches:
saddle-node bifurcations are available}
 \label{Ssn}

 We consider here some aspects of global continuation of the
$n$-branches of eigenvalues $\{\Lambda_l(n)\}$ obtained above.
Indeed, concerning the characteristic equation \ef{A7}, we do actually need
any advanced mathematical branching theory to see that
$\Lambda_l(n)$ emanate at $n=0$ from linear ones $\l_l$, since the
dependence on $n$ in \ef{A7} is analytic (linear), and must remain
like that for all $n>0$.

\ssk

\subsection{Saddle-node bifurcations for larger $n>0$ and $l \ge
2$}

\begin{figure}
\centering
\includegraphics[scale=0.75]{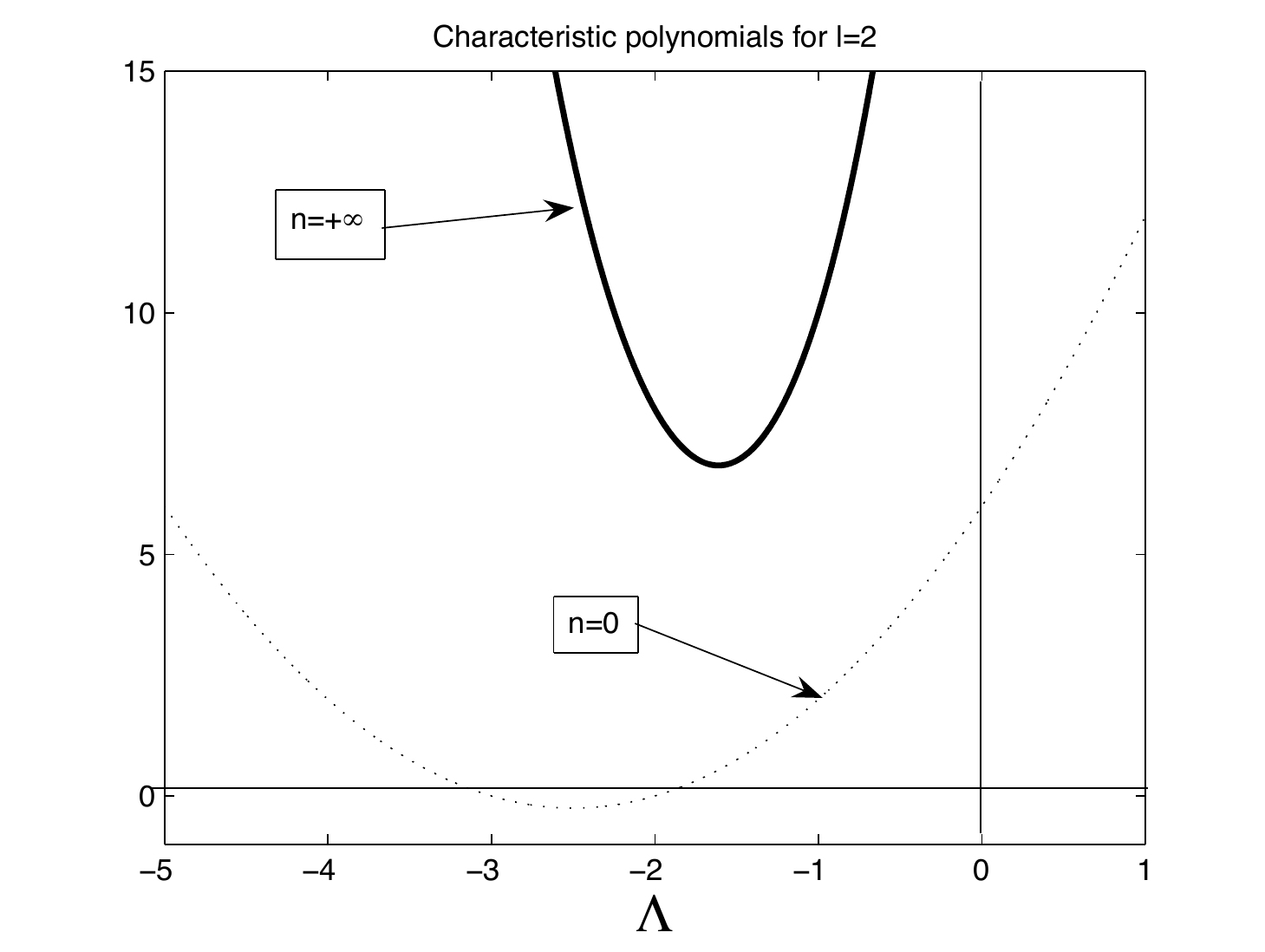}  
\vskip -.3cm
  \caption{Illustrations to characteristic polynomials \ef{A7} for $n=+\iy$ and $n=0$; $l=2$.}
 \label{Figl2}
\end{figure}

 To see whether
those eigenvalue branches persist for larger $n>0$, we first
consider the limit case $n=+\iy$, i.e., keeping two $O(n)$-terms
in \ef{A7}, we arrive at the asymptotic polynomial equation
 \be
 \label{p55}
F_l(\Lambda) \equiv [\Lambda^2+(2l+1)\Lambda+l(l+1)](\Lambda+l)^2+l^3(l-1)+2l(\Lambda+l)(\Lambda
l+l(l-1))=0.
 \ee
Next, to get a possibility of a saddle-node bifurcation in $n$,
let us restrict it to the simplest non-trivial case $l=2$, where
 \be
 \label{p56}
 F_2(\Lambda)=\Lambda^4+7\Lambda^3+26 \Lambda^2+46 \Lambda+36> 0 \quad \mbox{for all} \quad
 \Lambda \in \re,
 \ee
 which, thus, is strictly positive in $\re$. Moreover, the lower bound
 is rather impressive and convincing:
 \be
 \label{p57}
 F_2(\Lambda) \ge 6.84 \inB \re,
  \ee
  obtained through numerics.  
This is illustrated by
 Figure \ref{Figl2}, where the graph of \ef{p56} is presented by a
 bold-face line. By a dash line therein, we draw the
 characteristic polynomial \ef{eig11} for $n=0$, which has two clear
 roots $-2$ and $-3$ (i.e., $-l$ and $-l-1$ with $l=2$).

It follows that, in view of \ef{p56}, the characteristic equation
for the eigenvalues $\Lambda_l$ \ef{A7} does not have a solution
for all $n>0$ sufficiently large. Therefore, by the continuity
(analyticity), two branches of eigenvalues, emanating at $n=0$,
must be destroyed at a {saddle-node bifurcation} at some
$n=n_2^*>0$, so that, for $n>n_2^*$, no  real eigenvalues (and
hence, real eigenfunctions) for $l=2$ exist. Imaginary nonlinear
eigenvalues, as we know, are not associated with zeros sets of
real solutions of the $p$-Laplacian (or other operators).

We expect such a nonlinear
bifurcation phenomenon to exist for other $l \ge 3$, though the
corresponding bifurcation points $n_*^l$ can essentially depend on
$l$. In other words, for larger $n>0$, due to nonlinear properties
of the $p$-Laplacian, some ``multiple zero" structures may be
destroyed, while others may continue to exist. This will be
supported below by a piece of convincing numerical evidence.


\subsection{Sharp estimates of $n^*_l$ for $l=2,3,4, 5, 10, 100$}

In order to get those saddle-node eigenvalue bifurcation points, we
study numerically the characteristic equation \ef{A71} for $n>0$.
In view of \ef{A91}, we are not interested in the non-zero case
$l=0$.

\begin{figure}
\centering
\includegraphics[scale=0.75]{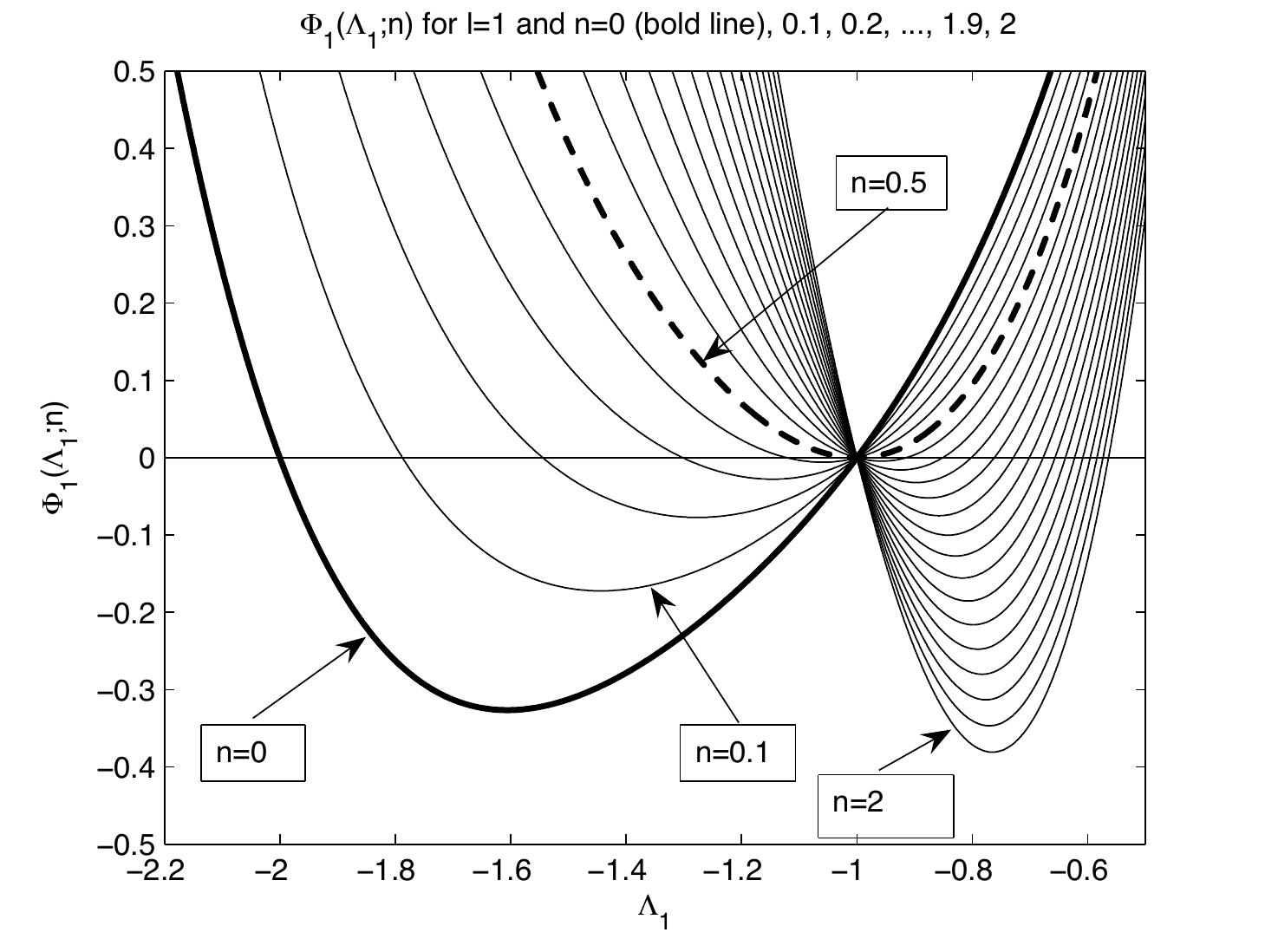}  
\vskip -.3cm
  \caption{Graph of the  characteristic polynomial \ef{A91} for $l=1$ and various values of $n=0, 0.1, ..., 1.9, 2$.}
 \label{FNew2}
\end{figure}

Again, by \ef{A91}, for $l=1$, we see that, for any $n>0$, there
exists 
$$\Lambda_1(n) =-1.$$ 
However, for any arbitrarily large
$n>0$, there exists another eigenvalue $\hat \Lambda_1(n)$ shown
in Figure \ref{FNew2}. It is interesting that, starting at the
bifurcation value 
$$\hat \Lambda_1(0)=-l-1=-2, \quad\hbox{precisely at}\quad \hat n^* = \frac 12,$$ 
the value of $\hat \Lambda_1(n)$ crosses -1 and
remains existent. Indeed, one can see from the characteristic
equation \ef{A71} that
 \be
 \label{bb21}
  \tex{
 l=1: \quad \mbox{there exists the double root} \,\,\,
 \hat \Lambda_1(n)= \Lambda_1(n)=-1, \,\,\, \mbox{iff} \,\,\, n= \frac 12.
 }
 \ee


\begin{figure}
\centering
\includegraphics[scale=0.75]{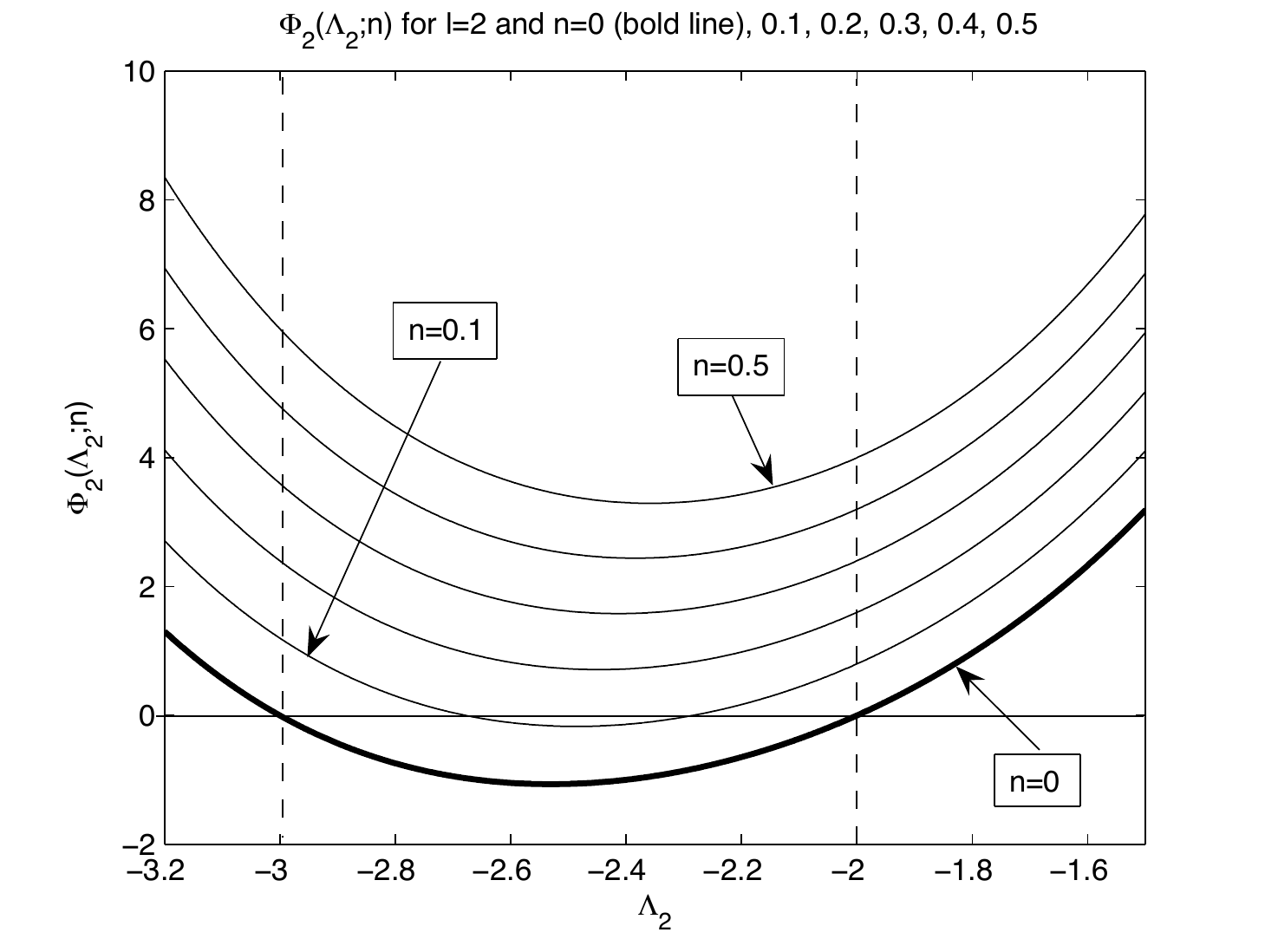}  
\vskip -.3cm
  \caption{Graph of the  characteristic polynomial \ef{A91} for $l=2$ and various values of $n=0, 0.1,0.2, 0.3,0.4,0.5$.}
 \label{FNew1}
\end{figure}

For $l=2$, the situation is different: Figure \ref{FNew1} clearly
shows existence of the saddle-node bifurcation at some $n_2^* \in
(0.1,0.2)$. More accurate numerics give the following estimates of
these bifurcation points and the corresponding coinciding
eigenvalues:
\be
\label{m21} 0.11912< n_2^* <0.11913 \whereA
\Lambda_2(n_2^*)=-2.4782...\, ,
 \ee
 so that, for $n> n^*_2$, there are no real eigenvalues of
 \ef{A91}.


\begin{figure}
\centering
\includegraphics[scale=0.75]{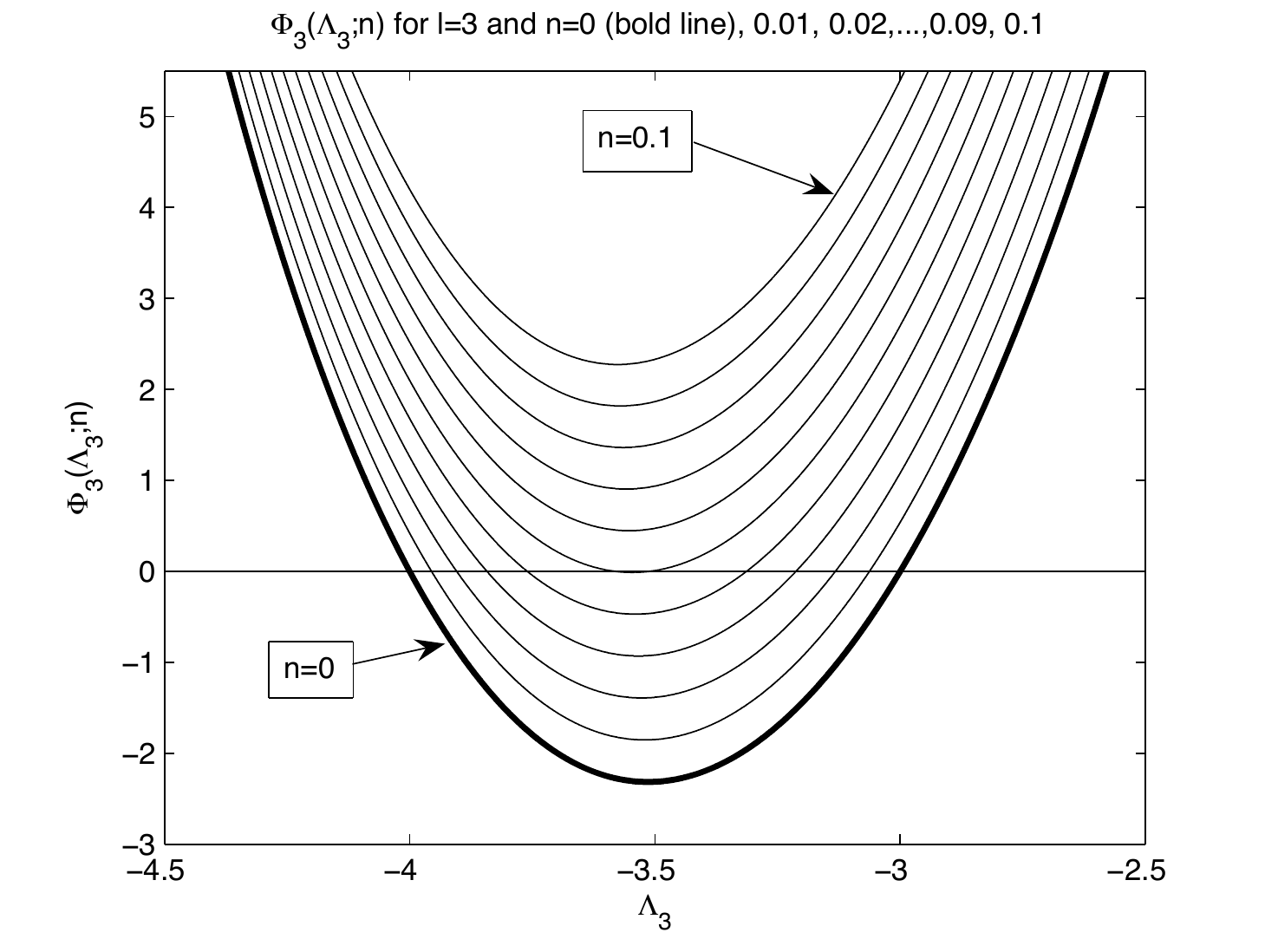}  
\vskip -.3cm
  \caption{Graph of the  characteristic polynomial \ef{A91} for $l=3$ and various values of $n=0, 0.01, ..., 0.1$.}
 \label{FNew4}
\end{figure}

\begin{figure}
\centering
\includegraphics[scale=0.75]{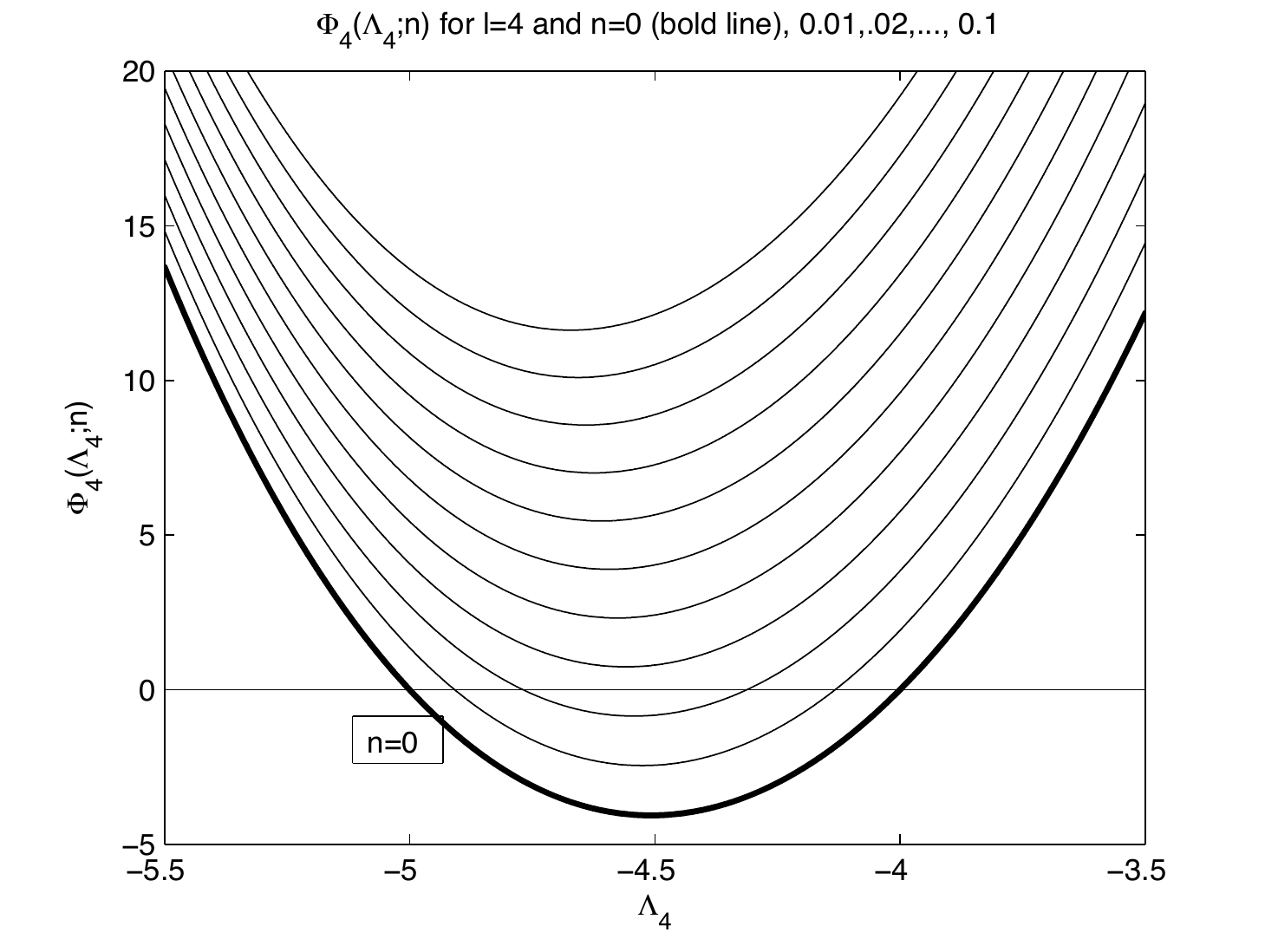}  
\vskip -.3cm
  \caption{Graph of the  characteristic polynomial \ef{A91} for $l=4$ and various values of $n=0, 0.001, ..., 0.01$.}
 \label{FNew3}
\end{figure}

Furthermore, Figures \ref{FNew4} and \ref{FNew3} show
similar bifurcation phenomena for $l=3$ and $l=4$ respectively.
The bifurcation points and eigenvalues satisfy:
 \be
 \label{m22}
 \begin{matrix}
 n_3^*=0.052292...,\,\,\,\Lambda_3(n_3^*)=-3.546... \ssk
  \\
   \andA
 n_4^*=0.025354..., \,\,\,  \Lambda_4(n_4^*)=-4.55514....\,
 .
 \end{matrix}
  \ee

 We expect that, for larger $l$, bifurcation values become even
smaller. For instance, our numerical estimates of the next
bifurcation point and coinciding (double) eigenvalues for $l=5$
are:
 \be
 \label{m23}
  n_5^*=0.014859 \whereA
 \Lambda_5(n^*_5)=-5.546...\,.
  \ee
Finally, we take larger $l=10$ and $l=100$:
 \be
 \label{m87}
 \begin{matrix}
 n_{10}^*=0.0036245... \, , \,\,\,
 \Lambda_{10}(n_{10}^*)=-10.5254... \ssk
  \\
  \andA
 n_{100}^*=0.00002555062... \, , \,\,\,
 \Lambda_{10}(n_{10}^*)=-100.5025... \, .
 \end{matrix}
  \ee
We see that, for $l=100$, the bifurcation-turning point of
eigenvalues $n$-branches occurs already at very small $n \sim
10^{-5}$, implying
 the
following final remarks naturally touching the cases $n>n_l^*$.



\subsection{\bf Problem of regularity at boundary points for nonlinear problems}

It should be mentioned that the problem of regularity at boundary points has been studied before 
using the so-called Wiener's  regularity test/criterion via concepts of potential-theoretic Bessel (Riesz) 
capacities to strongly elliptic equations. In that sense the problem is formulated in terms of a divergence series of capacities, 
measuring the thickness of the complement of the domain near the point of the boundary at which the regularity is analysed (in this case the origin 0); 
see \cite{Vic} for en extensive history of this problem applied 
to elliptic and parabolic equations, and the Kozlov--Maz'ya--Rossmann's monographs \cite{KMR1,KMR2} for 
just elliptic problems.

Therefore, following those arguments the principal extension of 
Wiener's-like capacity regularity test to a nonlinear degenerate $p$-Laplacian operator, with $p\in (1,N]$ 
was due to Maz'ya in 1970 \cite{Maz}, who extended it later on to have a sufficient capacity regularity condition to be optimal for any $p>1$.  

Note that the analysis of this kind of problem was completed for second order elliptic and parabolic equations by 
Wiener (1924) for the Dirichlet problem in \cite{Wien} and Petrovskii (1934), and extended in the 1960's-1970's for $2m$th-order PDEs by Kondriat'ev 
\cite{Kon2} and Maz'ya \cite{Maz,Maz2}. 
The key question was always to determine the optimal, and as sharp as possible, conditions on the ``shape" 
of the continuous  boundary under which the solutions is continuous at the boundary points.

Thus, one of the approaches was to treat the problem at the singular boundary point by \emph{blow-up evolution}
via approaching this ``singular" boundary point. Then, the rescaled pencil operators seem to be key. 
In that respect even though Konfriat'ev's paper was not devoted to regularity issues, it represented a novel idea
and involved the use of spectral properties of pencil operators, such as those we have used here, for describing 
the asymptotics of solutions near singularity boundary points assuming a kind of ``elliptic evolution" approach for elliptic problems.


\subsection{\bf Final comment on zeros for the $p$-Laplacian}

 In
this connection, it is important to mention that, for nonlinear
problems, ``blow-up" structures (at present, ``blow-up"  zero
structures) may be of different forms, especially, for the
$p$-Laplace operators. For instance, in \cite{BuGa}, it was shown
that blow-up in a quasilinear {\em Frank--Kamenetskii equation}
from combustion theory (a solid fuel model for $n=0$),
 \be
 \label{FK1}
 u_t=(|u_x|^n u_x)_x + \eee^u, \quad n >0,
  \ee
there exists an infinite sequence of critical exponents 
$$\{n_k\}  \to n_\iy<\iy,$$ 
so that the overall number of blow-up
 self-similar solutions changes by 1 when $n$ crosses each
 critical $n_k$. Recall that, for $n=0$, i.e., for the classic
 Frank--Kamenetskii equation from the 1930s, no such self-similar
 solutions exist at all, i.e., the above is a nonlinear
 phenomenon associated with the $p$-Laplace diffusion.

 Let us add that this phenomenon is associated with the increase
 of rotations on a phase plane admitted by ODEs for blow-up
 profiles, so that, for 
 $$n> n_\iy,$$ 
 there exist very complicated
 families of solutions, including those having infinite number of
 oscillations about a constant equilibrium; see also \cite{GPos}.
In other words, e.g., for sufficiently small 
$$0 < n <n_1,$$ 
there
exists a single ``nonlinear" blow-up structure, while all others
(a countable family) can be obtained via linearization and
matching, i.e., these  are closer to similar blow-up behaviour for
$n=0$. For 
$$n_2<n <n_3,$$ 
there are two self-similar profiles (the
rest are linearized), etc.

Similarly, in our nodal set case \ef{A40} completely coincides
with the linear representation, but until some maximal critical
value $n_*^l$. We do not know, whether, for $n>n_*^l$ there exist
other types of ``nonlinear multiple zeros" at $\t = + \iy$ for the
evolution equation \ef{A3}. Though, bearing in mind the {\em
finite propagation} property of the $p$-Laplacian in the evolution
sense and its infinitely oscillatory properties from
\cite{BuGa,GPos} for large exponents (also, a 2D problem, in the
rescaled variables $\{z,\t\}$, with $\t=-\ln(T-t)$, cf. $\t$ in
\ef{resc1}, similar to \ef{A3}), we can expect some complicated
{\em non-self-similar} zero structures of nodal sets that makes
the evolution completeness problem of multiple zeros extremely
difficult. Nevertheless, an ``accidental" 1-homogeneity (existence of an
``linear" invariant scaling) of \ef{A3} could make it slightly
easier.

 \section{Towards $p$-bi-Laplace equation}

Consider briefly \ef{d2}. Using the same rescaled variables
\ef{resc1} and bearing in mind the corresponding
$(z,\t)$-representation \ef{www.20} of the Laplacian, we arrive at
the equation
 \be
 \label{kk.1}
 \D_{(z,\t)}(|\D_{(z,\t)}w|^n \D_{(z,\t)} w)=0.
  \ee
  Next, performing all the technical computations, using the homogeneity of the final elliptic equation
  (recall the symmetry $w \mapsto Cw$, $C \in \re$, of \ef{kk.1}), we are looking
for solutions via nonlinear  eigenfunctions \ef{A40}, for which
 $$
 D_\t= \Lambda I,
 $$
 to get the corresponding nonlinear eigenvalue problem. The
 latter, for $n=0$, reduces to that for the bi-Laplacian one
 \cite{CGbiLap}, which is shown to admit four families of
 harmonic polynomials as eigenfunctions of a quadratic pencil of
 non self-adjoin operators. This allows:

 (i) to derive a nonlinear (polynomial) characteristic equation using the same
 assumption on the polynomial growth at infinity (associated with the
 analyticity of solutions of this ODE away from degeneracy points);

 (ii) to prove existence of nonlinear eigenfunctions for proper
 values of $\Lambda$;

 (iii) to establish a proper branching of nonlinear eigenfunctions at $n=0$
 from harmonic polynomials;

 (iv) to study, both analytically and numerically, those saddle-node
 $n$-bifurcation points for eigenvalues;

 \noi etc.

 In other words, these are quite similar to what we have
 done above for the $p$-Laplace equation, with more technical
 features, indeed.




\end{document}